\documentclass[11pt,a4paper]{article}
\usepackage{geometry}
\usepackage{graphicx}
\usepackage{bm,array}
\usepackage{comment}
\usepackage{caption}
\usepackage{subcaption}
\usepackage{tablefootnote}
\usepackage{makecell}
\usepackage{booktabs}
\usepackage{multirow}
\usepackage{float}
\usepackage{lscape}
\usepackage{graphicx}
\usepackage{multirow}
\usepackage{hhline}
\usepackage{graphicx}
\usepackage{bm,array}
\usepackage{comment}
\usepackage{caption}
\usepackage{subcaption}
\usepackage{makecell}
\usepackage{multirow}
\usepackage{float}
\usepackage{lscape}
\usepackage{svg}
\usepackage{framed}
\usepackage{apacite}
\usepackage{appendix}
\usepackage[ruled,vlined,linesnumbered]{algorithm2e}
\usepackage{algorithmic}
\usepackage[normalem]{ulem}
\usepackage{amssymb}
\usepackage{amsmath}
\usepackage{amsthm}

\usepackage[pdfpagemode={UseOutlines},bookmarks=true,bookmarksopen=true,
bookmarksopenlevel=0,bookmarksnumbered=true,hypertexnames=false,
colorlinks,linkcolor={blue},citecolor={blue},urlcolor={red},
pdfstartview={FitV},unicode,breaklinks=true]{hyperref}
\hypersetup{urlcolor=blue, colorlinks=true}

\usepackage{setspace}
\usepackage{vmargin}
\setmarginsrb{ 1.0in}  
{ 0.6in}  
{ 1.0in}  
{ 0.8in}  
{  20pt}  
{0.25in}  
{9pt}  
{ 0.3in}  

\usepackage{authblk}
\usepackage{bbm}
\usepackage{amsmath}
\usepackage{amsfonts}
\usepackage{amssymb}
\usepackage[round, sort,comma,authoryear]{natbib}
\usepackage{mdframed}

\usepackage{listings}
\newtheorem{theorem}{Theorem}
\newtheorem{assumption}{Assumption}

\newtheorem{lemma}{Lemma}

\newtheorem{proposition}{Proposition}

\usepackage{multirow}
\usepackage{hhline}
\usepackage{footnote}
\usepackage[ruled,vlined]{algorithm2e}
\usepackage{algorithmic}

\newcommand{\cP}{{\mathcal{P}}}

\newcommand{\cX}{{\mathcal{X}}}
\newcommand{\cF}{{\mathcal{F}}}

\newcommand{\cQ}{{\mathcal{Q}}}

\newcommand{\cO}{{\mathcal{O}}}

\newcommand{\cY}{{\mathcal{Y}}}
\newcommand{\cZ}{{\mathcal{Z}}}

\newcommand{\bB}{\textbf{B}}
\newcommand{\bC}{\textbf{C}}

\newcommand{\bx}{\textbf{x}}
\newcommand{\by}{\textbf{y}}

\newcommand{\bd}{\textbf{d}}
\newcommand{\bq}{\textbf{q}}

\newcommand{\bA}{\textbf{A}}

\newcommand{\bp}{\textbf{p}}
\newcommand{\ba}{\textbf{a}}
\newcommand{\bz}{\textbf{z}}

\newcommand{\bo}{\textbf{o}}
\newcommand{\bn}{\textbf{n}}

\newcommand{\btheta}{{\pmb{\theta}}}

\newcommand{\bbR}{\mathbb{R}}

\usepackage{mathtools}

\DeclarePairedDelimiter\floor{\lfloor}{\rfloor}

\usepackage{mathtools}

\newif\ifnotes\notestrue
\def\htien#1{}

\usepackage{xcolor}

\gdef\AQ#1{}
\gdef\CQ#1{}

\begin{document}

\newcolumntype{C}{>{\centering\arraybackslash}p{4em}}

\title{\textbf{Joint Binary-Continuous Fractional Programming: Solution Methods and Applications}}
\author[1]{Hoang Giang Pham}
\author[2]{Ngan Ha Duong}
\author[1]{Tien Mai}
\author[3]{Thuy Anh Ta}
\author[4]{Minh Hoàng Hà}
\affil[1]{\it\small
School of Computing and Information Systems, Singapore Management University, Singapore}
\affil[2]{\it\small
Department of Engineering, Computer Science and Mathematics, University of L'Aquila, Italia}
\affil[3]{\it\small
ORLab, School of Computing, Phenikaa University, Vietnam}
\affil[4]{\it\small
SLSCM and CADA, Faculty of Data Science and Artificial Intelligence, College of Technology, National Economics University, Vietnam}
\date{}

\maketitle

\begin{abstract}
In this paper, we investigate a class of non-convex sum-of-ratios programs relevant to decision-making in key areas such as product assortment and pricing, and facility location and cost planning. These optimization problems, characterized by both continuous and binary decision variables, are highly non-convex and challenging to solve. To the best of our knowledge, no existing methods can {efficiently} solve these problems to near-optimality with arbitrary precision. 
To address this challenge, we propose an innovative approach based on \textit{logarithmic transformations} and \textit{piecewise linear approximation (PWLA)} to approximate the nonlinear fractional program as a mixed-integer convex program with arbitrary precision, which can be efficiently solved using cutting plane (CP) or Branch-and-Cut (B\&C) procedures. Our method offers several advantages: it allows for a shared set of binary variables to approximate nonlinear terms and employs an optimal set of breakpoints to approximate other non-convex terms in the reformulation, resulting in an approximate model that is \textit{minimal in size}. Furthermore, we provide a theoretical analysis of the approximation errors associated with the solutions derived from the approximated problem. We demonstrate the applicability of our approach to constrained competitive joint facility location and cost optimization, as well as constrained product assortment and pricing problems. Extensive experiments on instances of varying sizes, comparing our method with several alternatives—including general-purpose solvers and more direct PWLA-based approximations—show that our approach consistently achieves superior performance across all baselines, particularly in large-scale instances.
\end{abstract}


\noindent\textbf{Notation:}
Boldface characters represent matrices (or vectors), and $a_i$ denotes the $i$-th element of vector $\ba$ if it is indexable. We use $[m]$, for any $m\in \mathbb{N}$, to denote the set $\{1,\ldots,m\}$. 

\section{Introduction}
We study the following non-convex optimization problem with binary and continuous variables
\begin{align}\label{prob:SoR-original}\tag{\sf SoR}
     \max_{\by,\bx}\left\{f(\by,\bx) = \sum_{t\in [T]}  \frac{a_t + \sum_{i\in [m]} y_i g^t_i(x_i) }{b_t + \sum_{i\in [m]} y_i h^t_i(x_i)}\Big|~~  (\bx,\by)\in \cZ\right\}
\end{align}
where $\by$ are binary and $\bx$ are continuous variables, $g^t_i(x),h^t_i(x)$ are univariate functions, i.e., $g^t_i(x),h^t_i(x) : \bbR\rightarrow \bbR$, $\forall t\in [T],i\in [m]$, noting that  $g^t_i(x),h^t_i(x)$, $t\in [T], i\in [m]$, are univariate functions and are not necessarily convex (or concave), and $\cZ$  is a feasible set of  $(\bx,\by)$ capturing some relations between the two sets of variables. Here, we assume that $\cZ$ can incorporate general linear constraints that capture business requirements on $\bx$ and $\by$, i.e.,
\[
\cZ = \left\{ (\by, \bx) \;\middle|\; x_i \in [l_i, u_i],\; y_i \in \{0,1\},\; \forall i \in [m],\; \text{and } \bA\by + \bB\bx \leq \bC \right\}.
\]
where $\bA\by + \bB\bx \leq \bC$ are some linear constraints on $\bx, \by$.
Such a sum-of-ratios problem arises from the use of discrete choice models \citep{McFa81,Trai03} to predict customer/adversary's behavior in decision-making and is known to be highly non-convex and challenging to solve, even when the binary variables $\by$ are fixed \citep{li2019product,duong2022joint}.  As far as we know, this is a first attempt to solve the aforementioned non-convex problems to near global optimality. The problem formulation above has several important applications in \textit{revenue  management and facility location,} as described below.

\noindent \textbf{Competitive facility location and cost optimization}. The formulation \eqref{prob:SoR-original} can be found along an active line of research on competitive maximum covering (or maximum capture) facility location problem with customers' random utilities \citep{benati2002maximum,FLO_Hasse2009MIP,mai2020multicut,lin2021branch,dam2022submodularity}. The problem refer to maximizing an expected customer demand, in a competitive market,  by locating new facilities and making decisions of the budget to spend on each opening facility,  assuming that customers make choice decisions according to a discrete choice model. When the costs are fixed, which is the focus of most of the works in the relevant literature,  researchers have shown that the facility location problem can be formulated as a mixed-integer linear program (MILP) \citep{benati2002maximum,freire2016branch, haase2014comparison}, or can be solved efficiently by outer-approximation algorithms \citep{mai2020multicut,Ljubic2018outer}. When the cost  optimization is considered
but the cost variables only take values from a discrete set, then it has been shown that the joint location and cost optimization problem can be converted to an equivalent facility location problem with binary variables and existing methods can apply \citep{qi2022sequential}. 
In contrast, if the cost variables are continuous,  the joint problem becomes highly non-convex and may have several local optima \citep{duong2022joint}. As far as we know, \cite{duong2022joint} is the only work to consider both facility location and cost optimization (with continuous cost variables). In this work, the authors state that the use of the standard mixed-logit model leads to a  intractable optimization problem with several local optimal solutions, and they instead propose to use  
a less-popular discrete choice framework, i.e., the multiplicative random utility maximization framework \citep{FosgBier09}.\textit{ So, the joint location and cost optimization  under the standard logit and mixed-logit model is still an open problem in the respective literature and we deal with it in this work.  }

\noindent \textbf{Product assortment and pricing optimization.} This problem refers to the problem of selecting a set of products and  making pricing decisions  to maximize an expected revenue, assuming that customers make  choice decisions according to a discrete choice model. Product assortment and pricing has been one of the most essential problems in revenue management and has received remarkable attention over the recent decades \citep{talluri2004revenue,vulcano2010om,rusmevichientong2014assortment,wang2018impact}. 
The joint assortment and price optimization problem under a (general) mixed-logit model (i.e., one of the most popular and general choice models in the literature) can be formulated in the form of \eqref{prob:SoR-original}. When the variables $\bx$ are fixed  and the objective function contains  only one ratio, the optimization problem can be solved in polynomial time under some simple settings, e.g. the problem is unconstrained or with a cardinality constraint \citep{rusmevichientong2010dynamic,talluri2004revenue}. 
When the objective function is a sum of ratios and the variable $\bx$ are fixed, the problem is generally NP-Complete even when there are only two fractions \citep{rusmevichientong2014assortment}. Approximate solutions, MILP and mixed-integer second order cone programming (SOCP) reformulations have been developed for this setting \citep{bront2009column,mendez2014branch,sen2018conic}. When only the pricing decisions are considered (i.e., the variables $\by$ are fixed)  and the objective function contains multiple ratios, the problem is highly non-convex and may have several local optima, with respect to both the prices and market shares \citep{li2019product}. Joint assortment and price optimization has been also studied in the literature \citep{wang2012capacitated}, but just under some simple settings (e.g., unconstrained on the prices and a cardinality constraint on the assortment, and the objective functions involves only one ratio). In general, as far as we know, \textit{in the context of assortment and price optimization, there is no global solution method to handle the joint problem with multiple ratios and general constraints. Our work is the first attempt to fill this literature gap.}


\noindent \textbf{Linear fractional programming. } Our work also relates to the literature of binary fractional programming and general fractional programming. In the context of  binary fractional programming, the problem is known to be NP-hard, even when there is only one ratio \citep{prokopyev2005complexity}. The problem is also hard to approximate \citep{prokopyev2005complexity}.    \cite{rusmevichientong2014assortment} show that for the unconstrained multi-ratio problem, there is no poly-time approximation algorithm that has an approximation factor better than $\cO(1/m^{1-\delta})$ for any $\delta >0$, where $m$ is the number of products. Exact solution methods for binary fractional programs include MILP reformulations \citep{mendez2014branch,haase2014comparison},  or Conic quadratic reformulations \citep{mehmanchi2019fractional,sen2018conic}. In fact, such MILP and Conic reformulations cannot be directly applied to our context due to the inclusion of  continuous variables.  Conversely, when the fractional program primarily deals with continuous variables (with fixed binary variables), it takes on a notably non-convex nature, leading to multiple local optima \citep{freund2001solving,gruzdeva2018solving}. Consequently, handling it exactly becomes  challenging. The general fractional program we are tackling, involving a combination of both binary and continuous variables, presents a particularly intricate problem to solve. \textit{To the best of our knowledge, there are currently no exact methods (except for some general-purpose solvers, which are typically inefficient) available in the respective literature  for achieving (near) optimal solutions in this context. }

\noindent \textbf{Piece-wise linear approximation (PWLA):}  Our work leverages a PWLA approach to simplify the objective function, leading to more tractable problem formulations. The literature on PWLA is extensive \citep{PWLA_lin2013review,PWLA_lundell2013refinement,PWLA_westerlund1998extended,PWLA_lundell2009some}, with various techniques integrated into state-of-the-art solvers for mixed-integer nonlinear programs \citep{GUROBI}. Although GUROBI's PWLA techniques offer methods to linearize certain types of nonlinear univariate functions, they are not directly applicable to solve the fractional program in \eqref{prob:SoR-original}. However, as discussed later, by reformulating problem~\eqref{prob:SoR-original} as a bilinear program, we demonstrate that GUROBI's PWLA capabilities can be applied. Nevertheless, this approach generally requires a large set of additional binary variables to approximate the nonlinear terms. In our experiments, we show that this method is consistently outperformed by our proposed solution techniques across most benchmark instances. It is also worth noting that PWLA techniques have been used to address MNL-based pricing problems~\citep{mai2023securing,Bose2022_NeurIPS}; however, these studies focus exclusively on single-ratio programs, rendering them unsuitable for the multi-ratio structure encountered in our setting.



\noindent \textbf{Our contributions. } We make the following contributions: 
\begin{itemize}
    \item[(i)] \textbf{Innovative approach based on log-transformation and PWLA.} We leverage a PWLA to tackle the challenging nonlinear fractional problem. While standard PWLA approaches typically require a large number of additional binary variables to approximate nonlinear terms—making them inefficient for large-scale problems—our goal is to develop a minimal-size approximation. To this end, we propose an innovative method that combines a logarithmic transformation with a sophisticated PWLA scheme to reformulate the original nonlinear fractional program into a mixed-integer convex program, which can be efficiently solved using Cutting Plan (CP) or Branch and Cut (B\&C) procedures. Our method offers several key advantages over direct PWLA-based methods:
    \begin{itemize}
        \item It allows for a shared set of binary variables to approximate all nonlinear terms $h^t_i(x_i)$ and $g^t_i(x_i)$, significantly reducing the model's complexity.
        \item It provides an optimal mechanism for selecting breakpoints to approximate some exponential terms in the formulation.
    \end{itemize}
    These features collectively yield a compact and scalable approximation model. Additionally, we explore several alternative (and more direct) PWLA-based approaches, including MILP- and SOCP-based reformulations and those utilizing GUROBI’s native PWLA functionality. We thoroughly discuss the comparative advantages of our approach (log-transformation + PWLA)  over these alternatives.

    \item[(ii)] \textbf{Theoretical guarantees.} We provide a theoretical analysis of the approximation guarantees offered by our approach. Specifically, we show that the combined use of log-transformation and PWLA yields an approximation error bounded by $\cO(\epsilon + 1/K)$, where $K$ denotes the number of breakpoints used to discretize the continuous variables $x_i$ for approximating the nonlinear terms $h^t_i(x_i)$ and $g^t_i(x_i)$, and $\epsilon$ represents the approximation error associated with exponential terms in the log-transformation. This result formalizes the intuition that increasing the granularity of the piecewise linear discretization improves solution quality, and it provides theoretical assurance that our approach can converge \textit{linearly} to an optimal solution of the original non-convex problem as $K$ increases and $\epsilon$ decreases.

    \item[(iii)] \textbf{State-of-the-art experimental performance.} We conduct comprehensive numerical experiments on instances of varying sizes, comparing our proposed method against multiple baselines, including a general mixed-integer nonlinear programming  solver and alternative PWLA-based techniques (e.g., MILP or  SOCP-based reformulations, and GUROBI's PWLA solver). The results clearly demonstrate the superiority of our approximation approach in producing near-optimal solutions for the original non-convex problem, consistently outperforming all baseline methods across the board.
\end{itemize}

In summary, we develop an innovative solution method based on logarithmic transformation and PWLA to obtain near-optimal solutions for the class of non-convex problems defined in~\eqref{prob:SoR-original}. Our approach not only provides solutions with provable approximation guarantees, but also offers a significantly more compact approximation model compared to alternative PWLA-based methods. Empirically, it achieves superior performance across all evaluated baselines in terms of both solution quality and computational efficiency. \textit{To the best of our knowledge, this work is the first to explore and develop global optimization techniques for several important classes of problems, including constrained product assortment and pricing, and constrained  facility location and cost planning.}

\noindent\textbf{Paper Outline.}
Section~\ref{sec:discr-approach} introduces our proposed approach based on logarithmic transformation and PWLA. Section~\ref{sec:performance-bounds} establishes theoretical performance guarantees for the approximation scheme. Section~\ref{sec:app} illustrates the applicability of our method to two important classes of problems: joint assortment and pricing optimization, and joint facility location and cost optimization. Section~\ref{sec:experiments} presents extensive numerical experiments to evaluate the effectiveness of our approach. Finally, Section~\ref{sec:concl} concludes the paper. Appendix \ref{appd:proofs} and \ref{appd:results} include technical proofs and additional experimental results. Supplementary Materials (Appendix C--D) present all baseline formulations.


\section{Solution Method: Log-transformation and PWLA}\label{sec:discr-approach}

To tackle the challenging mixed-integer nonlinear problem, our innovative solution method employs a \textbf{log-transformation approach} to simplify the fractional structure. We then utilize \textbf{PWLA} to linearize nonlinear terms and convexify the non-convex objective function. These steps enable us to approximate the original {binary-continuous non-convex program}  by  a \textbf{mixed-integer convex program (MICP)} with {arbitrarily high precision}. This reformulation allows the problem to be efficiently solved using CP or B\&C. In the following, we describe our approximation method step by step.

To begin, let us introduce the following mild assumptions which generally holds in all the aforementioned applications.  
\begin{assumption}\label{assm:as1}
The following assumptions hold:
\begin{itemize}
    \item [(i)] $b_t + \sum_{i\in [m]} y_i h^t_i(x_i)>0$ for all $\by\in \cY$, $\bx\in \cX$.
    \item[(ii)] $g^t_i(x_i),~h^t_i(x_i)$ are bounded, $\forall t\in [T], i\in [m]$.
    \item [(iii)] $h^t_i(x_i)$  and $g^t_i(x_i)$ are Lipschitz continuous, i.e., there exist  $L^{gt}_i, L^{ht}_i>0$ such that:
    \[|h^t_i(x_1) - h^t_i(x_2)| \leq L^{ht}_i |x_1-x_2| \text{ and } |g^t_i(x_1) - g^t_i(x_2)| \leq L^{gt}_i |x_1-x_2|,\forall x_1,x_2 \in [l_i,u_i].\]
\end{itemize}
\end{assumption}

\subsection{Converting to a Minimization Program} 
To convexify the objective function, we first reformulate it as a minimization problem. Specifically, we express the objective function as follows:
\begin{align}
f(\by,\bx) &= \sum_{t\in [T]}  \frac{a_t + \sum_{i\in [m]} y_i g^t_i(x_i) }{b_t + \sum_{i\in [m]} y_i h^t_i(x_i)} =  T\alpha -  \sum_{t\in [T]}  \frac{(\alpha b_t - a_t) + \sum_{i\in [m]} y_i (\alpha h^t_i(x_i) - g^t_i(x_i) )   }{b_t + \sum_{i\in [m]} y_i h^t_i(x_i)} \nonumber,
\end{align}
where $\alpha > 0$ is chosen to be sufficiently large such that:
\[
 \alpha(b_t + \sum_{i\in [m]} y_i h^t_i(x_i)) > a_t + \sum_{i\in [m]} y_i g^t_i(x_i).
\]
This choice of $\alpha$ is always feasible since the denominator $b_t + \sum_{i\in [m]} y_i h^t_i(x_i)$ remains positive, and the numerator $a_t + \sum_{i\in [m]} y_i g^t_i(x_i)$ is bounded from above (as per Assumption \ref{assm:as1}).

For notational simplicity, let us define \( u^t_i(x_i) = \alpha h^t_i(x_i) - g^t_i(x_i) \) and \( c_t = \alpha b_t - a_t \). Using these definitions, we can rewrite the objective function as:
\[
f(\by,\bx)  = T\alpha - \sum_{t\in [T]}  \frac{c_t + \sum_{i\in [m]} y_i u^t_i(x_i) }{b_t + \sum_{i\in [m]} y_i h^t_i(x_i)}.
\]
Consequently, we can reformulate the problem \eqref{prob:SoR-original} into the following minimization form:
\begin{align}
     \min_{\by\in \cY,\bx\in \cX }\left\{\cF(\by,\bx) = \sum_{t\in [T]}  \frac{c_t + \sum_{i\in [m]} y_i u^t_i(x_i) }{b_t + \sum_{i\in [m]} y_i h^t_i(x_i)} \Big|~~  (\bx, \by)\in \cZ\right\} \nonumber.
\end{align}
The motivation behind this reformulation is that, in the following steps, we will apply a logarithmic transformation to convert the sum-of-ratios program into a nonlinear program involving exponential and logarithmic terms. Under this minimization formulation, certain terms will become convex, facilitating the optimization process.

\subsection{Log-Transformation}
To simplify the fractional structure and convexify the objective function, we introduce the following logarithmic variables: $    n_t = \log \left( c_t + \sum_{i\in [m]} y_i u^t_i(x_i) \right) \text{ and }
    d_t = \log \left( b_t + \sum_{i\in [m]} y_i h^t_i(x_i) \right)$.
These transformations are always valid since both the numerator and denominator are strictly positive, i.e., $c_t + \sum_{i\in [m]} y_i u^t_i(x_i) > 0 \text{ and } b_t + \sum_{i\in [m]} y_i h^t_i(x_i) > 0, \quad \forall (\bx, \by) \in \cZ$.

Using these transformations, we reformulate the problem as:
\begin{align}
    \min_{\bx, \bn, \bd} \quad & \sum_{t\in [T]} e^{n_t - d_t} \nonumber\\
    \text{s.t.} \quad & e^{n_t} = c_t + \sum_{i\in [m]} y_i u^t_i(x_i), \quad \forall t\in[T], \nonumber\\
    & e^{d_t} = b_t + \sum_{i\in [m]} y_i h^t_i(x_i), \quad \forall t\in[T]. \nonumber
\end{align}
Since the objective function involves minimizing the terms $e^{n_t - d_t}$, we observe that $n_t$ should be maximized, while $d_t$ should be minimized as much as possible. Consequently, the equality constraints can be converted into inequalities: $e^{d_t} \leq c_t + \sum_{i\in [m]} y_i u^t_i(x_i)\text{ and } e^{n_t} \geq b_t + \sum_{i\in [m]} y_i h^t_i(x_i)$. Thus, we rewrite the problem as:
\begin{align}
    \min_{\bx, \bn, \bd} \quad & \sum_{t\in [T]} e^{n_t - d_t} \label{prob:exp-1}\tag{\sf LT1}\\
    \text{s.t.} \quad & e^{n_t} \geq c_t + \sum_{i\in [m]} y_i u^t_i(x_i), \quad \forall t\in[T], \label{ctr:e_nt}\\
    & e^{d_t} \leq b_t + \sum_{i\in [m]} y_i h^t_i(x_i), \quad \forall t\in[T]. \label{ctr:e_dt}
\end{align}
The above problem contains non-convex terms such as $u^t_i(x_i)$, $h^t_i(x_i)$, and the constraint \eqref{ctr:e_nt}, which require further convexification. In the following, we describe our discretization approach to approximately convexify the non-convex problem.

\subsection{Linearizing $u^t_i(x_i),~ h^t_i(x_i)$ via PWLA}\label{sec:pwla_u_h}
To convexify the nonlinear, nonconvex program in \eqref{ctr:e_nt} and \eqref{ctr:e_dt}, we employ PWLA to linearize the nonlinear terms \( u^t_i(x_i) \) and \( h^t_i(x_i) \). Typically, PWLA can be directly applied to linearize each univariate term \( u^t_i(x_i) \) or \( h^t_i(x_i) \) by representing it as a linear function of a set of additional binary and continuous variables, with separate sets of auxiliary variables introduced for different nonlinear terms. This approach is also implicitly implemented in state-of-the-art solvers with PWLA, such as GUROBI \citep{GUROBI}.

However, in our context, this standard approach requires introducing multiple additional binary variables—proportional to the number of nonlinear terms—which significantly increases the computational complexity of the approximation formulation. Our approach differs by using a\textbf{\textit{ shared set of binary variables}} for all univariate nonlinear terms \( u^t_i(x_i) \) and \( h^t_i(x_i) \), ensuring that the number of additional variables scales only with the number of original variables \( x_i \), rather than the number of nonlinear terms. This significantly reduces the computational burden while maintaining the accuracy of the approximation.

To describe the general idea,  we first let \( g(x): \mathbb{R} \to \mathbb{R} \) be a univariate function. Suppose \( g(x) \) is Lipschitz continuous over the interval \([l,u]\) with a Lipschitz constant \( L > 0 \), meaning that for all \( x_1, x_2 \in [l,u] \),
$|g(x_1) - g(x_2)| \leq L |x_1 - x_2|.$
For any \( K \in \mathbb{N} \), we discretize the interval \([l,u]\) into \( K \) equal subintervals of length \( \Delta = \frac{u-l}{K} \) and approximate \( x \) as: $
x \approx \widehat{x} = l + \Delta \left\lfloor \frac{x - l}{\Delta} \right\rfloor$.

To incorporate this approximation into a mixed-integer nonlinear programming (MINLP) formulation, we introduce binary variables \( z_k \in \{0,1\} \) for \( k \in [K] \) and approximate \( x \) as: $x \approx \widehat{x} = l + \Delta \sum_{k\in [K]} z_k$, where the binary variables satisfy the constraint \( z_k \geq z_{k+1} \), ensuring a unique active index. Specifically, for \( k^* = \lfloor (x - l) / \Delta \rfloor \), we enforce \( z_{k^*} = 1 \) and \( z_{k^*+1} = 0 \), effectively selecting the appropriate discrete approximation. Using this approximation, we can represent \( g(x) \) as a discrete linear function: $g(x) \approx g(\widehat{x}) = g(l) + \Delta \sum_{k\in [K]} \gamma^{g}_{k} z_k$, where \( \gamma^{g}_{k} \) represents the slope of \( g(x) \) in the interval \( [l + (k-1)\Delta, l + k\Delta] \), defined as:  
\[
\gamma^{g}_{k} = \frac{g(l + k\Delta) - g(l + (k-1)\Delta)}{\Delta}, \quad \forall k \in [K].
\]
By leveraging the Lipschitz continuity of \( g(x) \), we can bound the approximation error as  
$|g(x) - g(\widehat{x})| \leq L \Delta.$
Thus, as \( K \) increases, both \( \widehat{x} \) and \( g(\widehat{x}) \) converge linearly to \( x \) and \( g(x) \), respectively.

We now show how to use the  technique above to linearize the nonlinear terms $c_t + \sum_{i\in [m]} y_i u^t_i(x_i)$ and $ b_t + \sum_{i\in [m]} y_i h^t_i(x_i)$.
 For ease of notation, let us first denote $\Delta_i = {(u_i-l_i)/}{K}$. 
Given any $\bx\in \cX$, we approximate $x_i$, $i\in[m]$, by $l+ \Delta_i\floor{K(x_i-l)/(u_i-l_i)}$ and approximate the functions $u^t_i(x_i)$ and $h^t_i(x_i)$  by binary variables $z_{ik}\in \{0,1 \}$, $\forall i\in [m], k\in [K]$ as 
\[
    u^{t}_i(x_i) \approx   \widehat{u}^{t}_i(x_i)  = u^{t}_i(l_i) + \Delta_i\sum_{k\in [K]}  \gamma^{ut}_{ik}{z_{ik}} ;~~h^{t}_i(x_i) \approx \widehat{h}^{t}_i(x_i)  = h^{t}_i(l_i) + \Delta_i\sum_{k\in [K]}  \gamma^{ht}_{ik}{z_{ik}},
\]
where $\gamma^{ut}_{jk}$ and $\gamma^{ht}_{jk}$ are the slopes of $g^t(x_i)$ and $h^t(x_i)$ in $[l_i+ (k-1)\Delta_i;l_i+k\Delta_i]$, $\forall k\in[K]$, i.e., 
\begin{align}
    \gamma^{ut}_{ik} = \frac{1}{\Delta_i}\left(u^{t}_i(l_i+k\Delta_i) -u^{t}_i(l_i+(k-1)\Delta_i) \right),\; \forall k\in [K], \nonumber \\
    \gamma^{ht}_{ik} =\frac{1}{\Delta_i}\left(h^{t}_i(l_i+k\Delta_i) -h^{t}_i(l_i+(k-1)\Delta_i) \right),\; \forall k\in [K].  \nonumber 
\end{align}
Then, each  $\bx\in \cX$ can be written as: $x_i = l_i +\Delta_i \sum_{k\in [K]}z_{ik} + r_i$, where $r_i\in [0,\Delta_i)$ is used to capture the gap between $x_i$ and the binary approximation $l_i +\Delta_i \sum_{k\in [K]}z_{ik}$. 
We now can approximate \eqref{prob:exp-1} as the following  problem:
\begin{align}
    \min_{\bx, \by, \bz, \bn, \bd} \quad & \sum_{t\in [T]} e^{n_t - d_t} \label{prob:exp-2}\tag{\sf LT2}\\
    \text{s.t.} \quad & e^{n_t} \geq c_t+ \sum_{i\in [m]}y_i u^t_i(l_i) +  \sum_{i\in [m]}\Delta_i\sum_{k\in [K]} \gamma^{ut}_{ik} y_iz_{ik}, \quad \forall t\in[T], \label{ctr:pwl_e_nt}\\
    & e^{d_t} \leq  b_t+ \sum_{i\in [m]}y_i h^t_i(l_i) +  \sum_{i\in [m]}{\Delta_i}\sum_{k\in [K]} \gamma^{ht}_{ik} y_iz_{ik}, \quad \forall t\in[T], \label{ctr:pwl_e_dt}\\
    &  z_{ik}\geq z_{i,k+1},\quad \forall k\in [K-1], i\in [m], \label{ctr:exp-z-1} \\
&  {x}_i = l_i+\Delta_i \sum_{k\in [K]} z_{ik} + r_i,\quad \forall i\in [m],   \label{ctr:exp-z-2} \\
&  r_i \in [0,\Delta_i),\quad \forall i\in [m],\label{ctr:exp-z-3} \\
     & (\bx,\by)\in\cZ,\;\bz \in \{0,1\}^{m\times K}.\;    \nonumber
\end{align}
Constraints \eqref{ctr:pwl_e_nt} and \eqref{ctr:pwl_e_dt} involve bilinear terms $y_iz_{ik}$ which can be linearized by  introducing some additional binary variables $s_{ik} = y_iz_{ik}$  and linear constraints $s_{ik}\leq z_{ik},\; s_{ik}\leq y_{i},\; s_{ik}\geq z_{ik}+y_i-1,\;\forall i\in [m], k\in [K]$. However, upon closer examination, it can be shown that, under certain assumptions that typically hold, these bilinear terms can be linearized without introducing additional variables. To facilitate this point, we first let  $\cZ(\by)$  be the feasible set of the original problem \eqref{prob:SoR-original} with fixed $\by\in \cY$, i.e., $\cZ(\by) = \{\bx|\; \bx \in  \cX;\; (\bx,\by)\in \cZ \}$. We first introduce the following assumption that is needed for the result.
\begin{assumption}\label{assum:a2}
For any $(\bx,\by)\in \cZ$ we have $(\bx',\by)\in \cZ$ for all $(\bx',\by)\in \cZ$ and $\bx'\leq \bx$. 
\end{assumption}
The above assumption is not restrictive in our applications of interest. For instance, in the context of joint assortment and price optimization, it is sufficient to assume that prices lie within predefined lower and upper bounds or that a weighted sum of prices (with non-negative weight parameters) does not exceed an upper limit. Similarly, in cost optimization for location planning, one typically requires that the total cost does not exceed a specified budget, i.e.,  
$\sum_{i\in [m]} x_i \leq C.$
Under such constraints, Assumption \ref{assum:a2} is indeed satisfied.

Under Assumption \ref{assum:a2}, we show, in Proposition \ref{prop:approx-simplified}, that Problem \eqref{prob:exp-2} can be simplified by 
replacing $y_iz_{ik}$ by only $z_{ik}$ and adding  constraints  $y_i\geq  z_{i1}$, which implies that if  $y_i = 0$ then $z_{ik} = 0$ for all $i\in [m], k\in [K]$.   
\begin{proposition}\label{prop:approx-simplified}
Suppose Assumptions \ref{assm:as1} and \ref{assum:a2} hold, we have that \eqref{prob:exp-2} is equivalent to the following mixed-integer program:
\begin{align}
    \min_{\bx, \by, \bz, \bn, \bd} \quad & \sum_{t\in [T]} e^{n_t - d_t} \label{prob: exp-3}\tag{\sf LT3}\\
    \text{s.t.} \quad & e^{n_t} \geq c_t+ \sum_{i\in [m]}y_i u^t_i(l_i) +  \sum_{i\in [m]}\Delta_i\sum_{k\in [K]} \gamma^{ut}_{ik} z_{ik}, \quad \forall t \in [T], \label{ctr:simply_e_nt}\\
    & e^{d_t} \leq  b_t+ \sum_{i\in [m]}y_i h^t_i(l_i) +  \sum_{i\in [m]}{\Delta_i}\sum_{k\in [K]} \gamma^{ht}_{ik} z_{ik}, \quad \forall t \in [T], \label{ctr:simply_e_dt}\\
  &\pmb{ y_i\geq z_{i1}},\quad \forall  i\in [m],\nonumber \\
  & \text{Constraints } \eqref{ctr:exp-z-1}- \eqref{ctr:exp-z-2}- \eqref{ctr:exp-z-3}, \nonumber\\
     & (\bx,\by)\in\cZ,\;\bz \in \{0,1\}^{m\times K}.\;    \nonumber
\end{align}
\end{proposition}

In \eqref{prob: exp-3}, only the constraints \eqref{ctr:simply_e_nt} remain non-convex. To address this, we further employ a PWLA to linearize the exponential terms \( e^{n_t} \), for all \( t \in [T] \). Our approach is described below.

\subsection{Linearizing the Exponential Terms $e^{n_t}$}\label{sec:e_nt}
We describe a method to linearize the exponential terms \( e^{n_t} \). While \( e^{n_t} \) can be linearized using the same approach as outlined earlier—by dividing a feasible interval of \( n_t \) into equal subintervals and approximating \( e^{n_t} \) with additional binary variables—we propose a more efficient method that minimizes the number of breakpoints. The key idea is to carefully select breakpoints one by one, ensuring that the approximation error, i.e., the gap between \( e^{n_t} \) and its PWLA, does not exceed a predefined threshold \( \epsilon \). 

It is important to note that while the approach described below provides an optimal way to select breakpoints—resulting in a smaller number of breakpoints (and thus fewer additional binary variables) compared to the uniform-selection methods used for the nonlinear terms \( u^t_i(x_i) \) and \( h^t_i(x_i) \)—it is not suitable for linearizing \( u^t_i(x_i) \) and \( h^t_i(x_i) \) due to the following reasons:  
\begin{itemize}
    \item[(i)] \textit{This method relies on the convexity of the function \( e^x \), which may not hold for \( u^t_i(x_i) \) and \( h^t_i(x_i) \) under our general settings}
    \item[(ii)] \textit{Our goal is to linearize all nonlinear terms \( u^t_i(x_i) \) and \( h^t_i(x_i) \) (for all \( t \in [T], i \in [m] \)) using a shared set of additional binary variables. The optimal breakpoint selection procedure described below is not well-suited for this approach, as it requires each nonlinear term to be approximated by a separate (and optimal) set of binary variables. Consequently, it does not support the sharing of binary variables across multiple nonlinear terms.
.}
\end{itemize}

\paragraph{Linearization via Non-uniform Breakpoints}
For notational simplification, we denote $\cP(e^x|U,L,\bp)$ as a PWLA of the function $e^x$ when $x \in [L,U]$, with $\bp$ represents a vector of breakpoints in $[L,U]$ to construct the PWLA.
In the following, we describe our general approach to obtain $\cP(e^x|U,L,\bp)$\footnote{For notational simplicity, the notation used to describe the PWLA function $\cP(e^x|U,L,\bp)$ is independent of the main problem formulation and does not share the same interpretation.}. To start, we partition the interval $[L,U]$ into smaller sub-intervals using $H+1$ breakpoints $\bp = (p_1, \dots, p_{H+1})$ such that: $L = p_1 < p_2 < \dots < p_{H+1} = U$. By introducing additional variables $w_h \in [0,1]$ and $v_h \in \{0,1\}$ for all $h \in [H]$, we can construct a piecewise linear function to approximate $e^x$ as follows:
\[
\begin{cases}
    \cP(e^x|U,L,\bp) = e^L + \sum_{h \in [H]} \delta_k (p_{h+1} - p_h) v_h, \\
    x = L + \sum_{h \in [H]} (p_{h+1} - p_h) v_h, \\
    v_h \geq v_{h+1}, & \forall h \in [H], \\
    w_h \geq v_h, & \forall h \in [H], \\
    w_{h+1} \leq v_h, & \forall h \in [H-1],
\end{cases}
\]
where: $\delta_h = \frac{e^{p_{h+1}} - e^{p_h}}{p_{h+1} - p_h}, \quad \forall h \in [H]$, is the slope of $e^x$ in the interval $[p_h, p_{h+1}]$.

\paragraph{Optimal Selection of Breakpoints.}
We now describe how to \textbf{optimally} select the breakpoints $\bp = (p_h, h =\{1,\ldots,H+1\})$. Given an accuracy level $\epsilon$, our objective is to minimize $H$ while ensuring that: $\max_{x \in [L, U]} |\cP(e^x|U,L,\bp) - e^x| \leq \epsilon$.

To achieve this, we maximize the size of each sub-interval while ensuring the approximation error remains within $\epsilon$. Specifically, starting from each breakpoint $p_h$, we determine the next breakpoint $p_{h+1}$ such that the interval size $p_{h+1} - p_h$ is maximized while satisfying the error constraint. The approximation gap within $[p_h, p_{h+1}]$ is given by:
\[
\phi(x)  = e^{p_h} + (x - p_h) \frac{e^{p_{h+1}} - e^{p_h}}{p_{h+1} - p_h} - e^x.
\]
The worst-case error over the interval $[p_h, p_{h+1}]$ is obtained by solving the following convex optimization problem:
\[
\max_{x \in [p_h, p_{h+1}]} \phi(x).
\]
By taking the derivative of $\phi(x)$ and setting it to zero, the maximum deviation occurs at:
\[
x^* = \ln\left(\frac{e^{p_{h+1}} - e^{p_h}}{p_{h+1} - p_h}\right).
\]
The maximum gap is then given by:
\[
\theta(p_{h+1}) = e^{p_h} + \left(\ln\left(\frac{e^{p_{h+1}} - e^{p_h}}{p_{h+1} - p_h}\right) - p_h - 1\right) \frac{e^{p_{h+1}} - e^{p_h}}{p_{h+1} - p_h}.
\]
To ensure the approximation remains within the given error bound, we determine $p_{h+1}$ as:
\[
p_{h+1} = \arg\max_{t > p_h} \left\{t ~|~ \theta(t) \leq \epsilon \right\}.
\]
Since $\theta(t)$ is monotonically increasing in $t$, we can  efficiently determine $p_{h+1}$ using the following  binary search algorithm.
\begin{mdframed}[linewidth=1pt, roundcorner=5pt, backgroundcolor=gray!10]
\textbf{Binary Search Algorithm:}
\begin{itemize}
    \item \textbf{Step 1:} Initialize the lower bound $l = p_h$, upper bound $u = U$, and tolerance $\tau > 0$.
    \item \textbf{Step 2:} If $\theta(u) \leq \epsilon$, set $p_{h+1} = U$ and terminate.
    \item \textbf{Step 3:} Compute $w = (u + l)/2$. If $\theta(w) \leq \epsilon$, set $l = w$; otherwise, set $u = w$.
    \item \textbf{Step 4:} If $|u - l| \leq \tau$, return $p_{h+1} = l$ and terminate; otherwise, repeat \textbf{Step 3}.
\end{itemize}
\end{mdframed}

The above binary search algorithm converges exponentially to the optimal solution of:
\[
\max_{t > p_k} \left\{ t ~|~ \theta(t) \leq \epsilon \right\}.
\]
It can be shown that after $\log(1/\tau)$ iterations, the algorithm finds a solution $\widetilde{t}$ such that $|\widetilde{t} - t^*| \leq \tau$, where $t^*$ is the optimal solution.  We describe our general approach optimal breakpoints for constructing the approximation $\cP(e^x|U,L,\bp)$.
\begin{mdframed}[linewidth=1pt, roundcorner=5pt, backgroundcolor=gray!10]
\textbf{Constructing the Breakpoints}:\\
The breakpoints are determined iteratively as follows:
\begin{itemize}
    \item Initialize $p_1 = L$.
    \item Use the binary search procedure to find the next breakpoint $p_{h+1}$.
    \item Terminate when $p_{h+1} = U$.
\end{itemize}
\end{mdframed}
Since the PWLA gap is optimized within each sub-interval, this method provides the optimal number of breakpoints. That is, no alternative set of breakpoints exists with a smaller $H$ while satisfying:
\[
\max_{x \in [L, U]} |\cP(e^x|U,L,\bp) - e^x| \leq \epsilon.
\]
In practice, choosing a very small threshold $\tau$ ensures near-optimality. Due to the exponential convergence of the binary search, the method terminates after only a few iterations, even for very small $\tau$.

We characterize some key properties of the PWLA function \( \mathcal{P}(e^x \mid U, L, \bp) \), which are important for understanding its behavior as well as for establishing performance guarantees when using this approximation in the overall nonlinear fractional program.

\begin{theorem}\label{theorem1}
The following properties hold:
\begin{itemize}
    \item[(i)] The PWLA function \( \mathcal{P}(e^x \mid U, L, \bp) \) is strictly monotonically increasing in \( x \) for all \( x \in [L,U] \). As a result, there exists a well-defined inverse function \( \mathcal{P}^{-1}(z, \bp) \) such that, for any \( z \in [e^L,e^U] \),
    \[
    \mathcal{P}(e^{\mathcal{P}^{-1}(z, \bp)} \mid U, L, \bp) = z.
    \]
    \item[(ii)] The number of breakpoints generated by the above procedure can be bounded as:
    \[
    H \leq \frac{e^U (U - L)}{\epsilon} + 1,
    \]
    implying that the breakpoint optimization procedure will terminate after at most $O(\frac{U - L}{\epsilon})$ iterations.
\end{itemize}
\end{theorem}

\subsection{Mixed-integer Convex Approximation}\label{sec:micp}
Combining all the techniques described above, we formulate a tractable approximation of the joint binary-continuous fractional program in \eqref{prob:SoR-original}. For notational simplification, let \( \cP(e^{n_t}|L_t, U_t,\bp^t) \) denote the PWLA of \( e^{n_t} \) for each \( t \in [T] \), where the breakpoints \( \bp^t \) are optimally constructed using the procedure described earlier. Using this approximation, we can approximate \eqref{prob:SoR-original} with the following mixed-integer convex program:
\begin{mdframed}[linewidth=1pt, roundcorner=5pt, backgroundcolor=gray!10]
    \begin{align}
    \min_{\bx, \by, \bz, \bn, \bd} \quad & \sum_{t\in [T]} e^{n_t - d_t} \label{prob:approx-convex-exp}\tag{\sf MICP1}\\
    \text{s.t.} \quad & \cP(e^{n_t}|L_t,U_t,\bp^t) \geq c_t+ \sum_{i\in [m]}y_i u^t_i(l_i) +  \sum_{i\in [m]}\Delta_i\sum_{k\in [K]} \gamma^{ut}_{ik} z_{ik}, \quad \forall t \in [T], \label{ctr:linear-P}\\
    & e^{d_t} \leq  b_t+ \sum_{i\in [m]}y_i h^t_i(l_i) +  \sum_{i\in [m]}{\Delta_i}\sum_{k\in [K]} \gamma^{ht}_{ik} z_{ik}, \quad \forall t \in [T], \label{ctr:e-dt}\\
  &{ y_i\geq z_{i1}},\quad \forall  i\in [m],\nonumber \\
  & \text{Constraints } \eqref{ctr:exp-z-1}- \eqref{ctr:exp-z-2}- \eqref{ctr:exp-z-3}, \nonumber\\
     & (\bx,\by)\in\cZ,\;\bz \in \{0,1\}^{m\times K},\;    \nonumber
\end{align}
\end{mdframed}
where $L_t$ and $U_t$ are lower bound and upper bound of $n_t$, which are also lower   and upper bounds of $\log(c_t + \sum_{i\in [m]} y_i u^t_i(x_i))$, which can be estimated  quickly. 

Note that \eqref{prob:approx-convex-exp} is only valid under Assumption \ref{assum:a2}. If this assumption does not hold, as discussed earlier, we can introduce additional variables \( s_{ik} \) to represent the terms \( y_i z_{ik} \) and linearize these bilinear terms using McCormick inequalities \citep{mccormick1976computability}.

The constraints in \eqref{ctr:linear-P} are linear since $\mathcal{P}(e^{n_t} \mid L_t, U_t, \bp^t)$ is a piecewise linear function. Therefore, the nonlinear program in \eqref{prob:approx-convex-exp} has a convex objective and convex constraints, which can generally be solved to optimality using CP or B\&C methods.   The general idea is to reformulate a master problem with linear constraints and solve it iteratively by adding valid cuts that approximate the convex constraints and objective. Specifically, we reformulate  \eqref{prob:approx-convex-exp} as follows: 
 \begin{align}
    \min_{\bx, \by, \bz, \bn, \bd, \btheta, \pmb{\psi}} \quad & \sum_{t\in [T]} \theta_t \label{prob:approx-convex-exp-1}\tag{\sf MICP2}\\
    \text{s.t.} \quad & \cP(e^{n_t}|L_t,U_t,\bp^t) \geq c_t+ \sum_{i\in [m]}y_i u^t_i(l_i) +  \sum_{i\in [m]}\Delta_i\sum_{k\in [K]} \gamma^{ut}_{ik} z_{ik}, \quad \forall t \in [T], \label{ctr:cP}\\
    &\psi_t \leq  b_t+ \sum_{i\in [m]}y_i h^t_i(l_i) +  \sum_{i\in [m]}{\Delta_i}\sum_{k\in [K]} \gamma^{ht}_{ik} z_{ik}, \quad \forall t \in [T],\label{ctr:psi-t-bt}\\
  &{ y_i\geq z_{i1}},\quad \forall  i\in [m],\label{ctr:y-z1} \\
  &\theta_t \geq e^{n_t-d_t} ,\quad \forall  t\in [T],~\label{ctr:theta-t}\\
  &\psi_t \geq e^{d_t} ,\quad \forall  t\in [T],~\label{ctr:psi- t}\\
  & \text{Constraints } \eqref{ctr:exp-z-1}- \eqref{ctr:exp-z-2}- \eqref{ctr:exp-z-3}. \nonumber
\end{align}
A master program can be defined by removing constraints \eqref{ctr:theta-t} and \eqref{ctr:psi- t} from \eqref{prob:approx-convex-exp-1}, and replacing them with a set of linear cuts that are added iteratively. Specifically, at each iteration, we solve the master problem to obtain a solution 
$(\overline{\bx}, \overline{\by}, \overline{\bz}, \overline{\bn}, \overline{\bd}, \overline{\btheta}, \overline{\pmb{\psi}}),$
and add the following gradient-based valid cuts to the master problem:
\begin{align}
\theta_t &\geq  e^{\overline{n}_t - \overline{d}_t} \left( 1 + (n_t - d_t) - (\overline{n}_t - \overline{d}_t) \right), \label{ctr-OA-1}\\
\psi_t &\geq  e^{\overline{d}_t} \left( 1 + d_t - \overline{d}_t \right). \label{ctr-OA-2}
\end{align}
It can be seen that the number of linear cuts added to the master problem is proportional to $T$. In addition, we can further enhance the CP or B\&C process by the following valid cuts. To present these valid cuts, let us denote:
\[
\eta_t(\by,\bz) = b_t + \sum_{i \in [m]} y_i h^t_i(l_i) + \sum_{i \in [m]} \Delta_i \sum_{k \in [K]} \gamma^{ht}_{ik} z_{ik}.
\]

\begin{proposition}\label{prop:SOCP-2}
Given any solution candidate 
$(\overline{\bx}, \overline{\by}, \overline{\bz}, \overline{\bn}, \overline{\bd}, \overline{\btheta}, \overline{\pmb{\psi}}),$
the following cuts are valid for \eqref{prob:approx-convex-exp-1}:
\begin{equation}
d_t \leq \log(\eta_t(\overline{\by}, \overline{\bz})) + \frac{
\sum_{i \in [m]} h^t_i(l_i) (y_i - \overline{y}_i)
}{\eta_t(\overline{\by}, \overline{\bz})} + \frac{
\sum_{i \in [m]} \sum_{k \in [K]} \Delta_i \gamma^{ht}_{ik} (z_{ik} - \overline{z}_{ik})
}{\eta_t(\overline{\by}, \overline{\bz})}, \quad \forall t \in [T].\label{ctr-OA-3}
\end{equation}
\end{proposition}

We describe our CP algorithm as follows. We first define the master problem derived from \eqref{prob:approx-convex-exp-1}, in which the nonlinear constraints are removed:
\begin{align}
\min_{\bx, \by, \bz, \bn, \bd, \btheta, \pmb{\psi}} \quad & \sum_{t \in [T]} \theta_t \label{prob:master-exp-1} \tag{\sf Master}\\
\text{s.t.} \quad & \text{Constraints } \eqref{ctr:exp-z-1}- \eqref{ctr:exp-z-2}- \eqref{ctr:exp-z-3}-\eqref{ctr:cP}-\eqref{ctr:psi-t-bt}-\eqref{ctr:y-z1}. \nonumber
\end{align}
The following procedure outlines the specific steps of the CP algorithm:

\begin{mdframed}[linewidth=1pt, roundcorner=5pt, backgroundcolor=gray!10]
\textbf{CP Procedure:}
\begin{itemize}
    \item \textbf{Step 1 (Find a solution candidate):} Solve the master problem \eqref{prob:master-exp-1} to obtain a solution candidate $(\overline{\bx}, \overline{\by}, \overline{\bz}, \overline{\bn}, \overline{\bd}, \overline{\btheta}, \overline{\pmb{\psi}})$.
    \item \textbf{Step 2 (Check feasibility):} Verify whether the solution candidate satisfies the mixed-integer nonlinear program \eqref{prob:approx-convex-exp}, i.e., check if
    \[
    \overline{\theta}_t \geq e^{\overline{n}_t - \overline{d}_t} - \xi \quad \text{and} \quad
    \overline{\psi}_t \geq e^{\overline{d}_t} - \xi,
    \]
    for all $t \in [T]$ and a given threshold $\xi > 0$ chosen as a stopping condition. If these inequalities hold, terminate the CP procedure and return $(\overline{\bx}, \overline{\by})$ as the optimal solution. Otherwise, proceed to \textbf{Step 3}.
    \item \textbf{Step 3 (Add cuts and iterate):} Add the linear cuts described in \eqref{ctr-OA-1}, \eqref{ctr-OA-2}, and \eqref{ctr-OA-3} to the master problem \eqref{prob:master-exp-1}, and return to \textbf{Step 1}.
\end{itemize}
\end{mdframed}

The CP method described above is guaranteed to solve the nonlinear program \eqref{prob:approx-convex-exp} to global optimality~\citep{OA_Duran1986outer,OA_Bonami2008algorithmic}. The valid cuts introduced, such as those in \eqref{ctr-OA-1}, \eqref{ctr-OA-2}, and \eqref{ctr-OA-3}, can also be incorporated into a B\&C procedure for solving the same problem. The core idea in the B\&C framework is similar: at each node of the branch-and-bound tree, a relaxed version of the nonlinear problem is considered (i.e., the master problem without nonlinear constraints). Valid cuts are then iteratively added to this relaxed master problem, thereby tightening the feasible region and refining the approximation of the nonlinear constraints. Solving this refined master problem provides an upper bound on the optimal objective value at that node~\citep{Ljubic2018outer}.  By systematically branching on integer variables and incorporating these valid cuts at each node, the B\&C method can efficiently navigate the solution space while maintaining valid bounds, ultimately converging to the global optimum of the original mixed-integer nonlinear program \eqref{prob:approx-convex-exp}.

\subsection{Alternative Tractable Approximations}\label{sec:milp_socp}
In the above, we presented our approach to solve the joint binary-continuous program by leveraging a log-transformation, PWLA using shared binary variables, and an optimal way to approximate exponential functions. The goal was to derive an approximation program that is optimal in size. In fact, PWLA can be applied to approximate and reformulate the joint problem in a more \textit{direct} and \textit{straightforward} manner. In the following, we discuss several alternative approaches that use PWLA to approximate the joint problem with a reformulated model that can be solved using existing solvers such as Gurobi.

\paragraph{MILP and SOCP Approximations via PWLA.}
We first note that PWLA enables the linearization of the nonlinear numerator and denominator in the original objective function \eqref{prob:SoR-original}. Consequently, one can approximate the original nonlinear fractional program via a linear fractional program, which can then be transformed into a MILP or second-order cone program (SOCP) using established techniques. Specifically, via the PWLA method described in Section~\ref{sec:pwla_u_h}, we can approximate \eqref{prob:SoR-original} by the following binary linear fractional program:
\begin{align}
    \min_{\bx, \by, \bz, \bn, \bd} \quad & \sum_{t \in [T]} \frac{ c_t + \sum_{i \in [m]} y_i u^t_i(l_i) + \sum_{i \in [m]} \Delta_i \sum_{k \in [K]} \gamma^{ut}_{ik} z_{ik} }{ b_t + \sum_{i \in [m]} y_i h^t_i(l_i) + \sum_{i \in [m]} \Delta_i \sum_{k \in [K]} \gamma^{ht}_{ik} z_{ik} } \label{prob:linear-frac-1} \tag{\sf LF}\\
    \text{s.t.} \quad & \text{Constraints } \eqref{ctr:exp-z-1}- \eqref{ctr:exp-z-2}- \eqref{ctr:exp-z-3}- \eqref{ctr:y-z1}, \nonumber \\
    & (\bx, \by) \in \mathcal{Z},\quad \bz \in \{0,1\}^{m \times K}. \nonumber
\end{align}
Note that the above formulation is only valid under Assumption \ref{assum:a2}. If this assumption does not hold, one can simply introduce additional variables to linearize the terms \( y_i z_{ik} \). However, for simplicity, we retain the formulation in \eqref{prob:linear-frac-1} as is.

Since \eqref{prob:linear-frac-1} is a binary linear fractional program, it can be conveniently reformulated as either a MILP or a mixed-integer SOCP using McCormick inequalities (see Supplementary Materials--Appendix C). A key distinction from prior work is that the coefficients in the denominator of \eqref{prob:linear-frac-1} can be negative when the function \( h^t_i(x_i) \) is decreasing in \( x_i \) (as in assortment and pricing problems). As a result, standard CONIC reformulations in~\cite{sen2018conic} do not apply. We explicitly discuss this issue in Supplementary Materials--Appendix C.2.

The main \textit{disadvantage} of the MILP and SOCP reformulations, compared to our log-transformation-based approximation in \eqref{prob:approx-convex-exp}, is that they require a large number of additional variables. The size of the reformulated model grows rapidly with \( m \) and \( K \). Moreover, the linearization of bilinear terms using McCormick inequalities leads to weak continuous relaxations.

\paragraph{Gurobi's PWLA.}
It is important to note that PWLA has been incorporated into several state-of-the-art solvers, such as Gurobi, to efficiently handle mixed-integer nonlinear programs~\citep{GUROBI}. While such PWLA-based techniques cannot be directly used to solve fractional programs like \eqref{prob:linear-frac-1}—since solvers such as Gurobi do not natively support fractional objectives—we can leverage their ability to handle bilinear terms by reformulating the fractional program as a bilinear one. In particular, Gurobi's PWLA can then be used to linearize the nonlinear functions \( g^t_i(x_i) \) and \( h^t_i(x_i) \).

Specifically, define the following quantities:
\begin{align*}
    o_t = c_t + \sum_{i \in [m]} y_i u^t_i(x_i);~~
    q_t = b_t + \sum_{i \in [m]} y_i h^t_i(x_i);~~
    \theta_t = \frac{n_t}{d_t},
\end{align*}
and rewrite the original nonlinear fractional program \eqref{prob:linear-frac-1} as the following bilinear program:
\begin{align}
    \min_{\bx, \by, \bz, \bo, \bq, \pmb{\theta}} \quad & \sum_{t \in [T]} \theta_t \label{prob:linear-frac-Blinear} \tag{\sf LFBL}\\
    \text{s.t.} \quad & o_t \geq c_t + \sum_{i \in [m]} y_i u^t_i(x_i), \quad \forall t \in [T],\nonumber \\
                      & q_t \leq b_t + \sum_{i \in [m]} y_i h^t_i(x_i), \quad \forall t \in [T],\nonumber \\
                      & \theta_t \cdot q_t = o_t, \quad \forall t \in [T],\nonumber \\
                      & \text{Constraints } \eqref{ctr:exp-z-1}- \eqref{ctr:exp-z-2}- \eqref{ctr:exp-z-3}- \eqref{ctr:y-z1}, \nonumber \\
                      & (\bx, \by) \in \mathcal{Z}, \quad \bz \in \{0,1\}^{m \times K}. \nonumber
\end{align}
We can now linearize the bilinear terms \( y_i u^t_i(x_i) \) using McCormick inequalities (or handle them directly by Gurobi) and further apply Gurobi’s built-in PWLA to approximate the univariate nonlinear functions \( u^t_i(x_i) \). The same technique is applied to the terms \( y_i h^t_i(x_i) \). By combining these linearizations, the bilinear program \eqref{prob:linear-frac-Blinear} becomes a mixed-integer linear approximation that can be solved to near-optimality using Gurobi.

The key difference between this approach and our PWLA method  is that Gurobi’s PWLA approximates each exponential function via separate sets of built-in piecewise linear constraints, rather than using a single set of variables \( \{z_{ik},~i\in[m], k\in [K] \} \) as we do. As a result, the number of additional binary variables in Gurobi’s PWLA is proportional to both the number of original variables \( x_i \) (for all \( i \in [m] \)) and the number of exponential terms in the objective function. In contrast, our PWLA approach only discretizes the original continuous variables \( \{x_i,~i \in [m]\} \), resulting in a more compact formulation and  fewer binary variables than the Gurobi-based PWLA.
    

In the experimental section, we will compare our log-transformation approach against all the aforementioned methods. Our results demonstrate that the log-transformation combined with PWLA consistently outperforms other approaches, especially for large-scale instances.

\section{Performance Guarantees}
\label{sec:performance-bounds}
In this section, we analyze the approximation errors yielded by solving the approximate program \eqref{prob:approx-convex-exp}. Our approximation scheme consists of \textit{two layers of approximation: }
\begin{itemize}
    \item[(i)] Approximating the univariate functions \( u^t_i(x_i) \) and \( h^t_i(x_i) \) using a common set of uniform breakpoints ($K+1$ breakpoints for each variable $x_i$, $i\in [m]$).
    \item[(ii)] Approximating the exponential function \( e^{n_t} \), for all \( t \in [T] \), by PWLA with an optimal set of breakpoints $\bp$, ensuring the guarantee: $\max_{x \in [L, U]} \left| \mathcal{P}(e^x \mid U, L, \bp) - e^x \right| \leq \epsilon$.
\end{itemize}
Thus, our goal is to establish an upper bound for the approximation errors as a function of the number of breakpoints \( K \) and the accuracy level \( \epsilon \). These analyses provide insights into how the parameter \( K \) (the number of pieces) and the accuracy level \( \epsilon \) affect the quality of the approximation, as well as guidelines for selecting an appropriate \( K \) and \( \epsilon \) to achieve the desired solution accuracy. This helps balance the trade-off between computational efficiency and approximation quality.


First, to facilitate our analysis and simplify notation, let us define:
\[
    u^t(\by,\bx) = c_t + \sum_{i\in [m]} y_i u^t_i(x_i),
    \quad \text{and} \quad
    h^t(\by,\bx) = b_t + \sum_{i\in [m]} y_i h^t_i(x_i).
\]
We denote their PWLAs based on our discretization technique as \( \widehat{u}^t(\by,\bx) \) and \( \widehat{h}^t(\by,\bx) \), given by:
\begin{align}
    \widehat{u}^t(\by,\bx) &=  c_t + \sum_{i\in [m]} y_i u^t_i(l_i) + \sum_{i\in [m]} \Delta_i \sum_{k\in [K]} \gamma^{ut}_{ik} z_{ik}, \\
    \widehat{h}^t(\by,\bx) &=  b_t + \sum_{i\in [m]} y_i h^t_i(l_i) + \sum_{i\in [m]} \Delta_i \sum_{k\in [K]} \gamma^{ht}_{ik} z_{ik}.
\end{align}
Here, the variables \( (\bx, \by, \bz) \) satisfy the constraints \eqref{ctr:exp-z-1}, \eqref{ctr:exp-z-2}, and \eqref{ctr:exp-z-3}, and  condition \( y_i \geq z_{i1} \) for all \( i \in [m] \).  Let \( \xi^u_t(\by,\bx) \) and \( \xi^h_t(\by,\bx) \) represent the approximation errors introduced by these PWLAs:
\[
\xi^u_t(\by,\bx) = u^t(\by,\bx) - \widehat{u}^t(\by,\bx), \quad \text{and} \quad
\xi^h_t(\by,\bx) = h^t(\by,\bx) - \widehat{h}^t(\by,\bx).
\]
We now establish bounds for these approximation errors.
\begin{lemma}\label{lem:1}
    For any \( (\by,\bx) \in \cZ \), the approximation errors satisfy:
    \[
        |\xi^u_t(\by,\bx)| \leq  \sum_{i\in [m]} (\alpha L^{ht}_i + L^{gt}_i) \frac{u_i - l_i}{K}, 
        \quad \text{and} \quad
        |\xi^h_t(\by,\bx)| \leq  \sum_{i\in [m]} L^{ht}_i \frac{u_i - l_i}{K}.
    \]
\end{lemma}

We further denote \( \xi_t(n_t) \) as the gap between \( e^{n_t} \) and its PWLA \( \mathcal{P}(e^{n_t} \mid L_t,U_t,\epsilon) \), i.e., 
\[
\xi_t(n_t) = e^{n_t}-\mathcal{P}(e^{n_t} \mid L_t,U_t,\epsilon).
\]
To establish a bound for the approximation error yielded by the approximation \eqref{prob:approx-convex-exp}, we rewrite the approximation program \eqref{prob:approx-convex-exp} equivalently as:
\begin{align}
    \min_{\bx, \by, \bn, \bd} \quad & \sum_{t\in [T]} e^{n_t - d_t} \label{prob:proof-1} \\
    \text{s.t.} \quad &  e^{n_t} - \xi^t(n_t) \geq u^t({\by,\bx}) - \xi^u_t({\by,\bx}), \quad \forall t \in [T], \nonumber\\
    & e^{d_t} \leq  h^t(\by, \bx) - \xi^h_t(\by,\bx), \quad \forall t \in [T], \nonumber\\
    & (\bx,\by) \in \cZ. \nonumber
\end{align}
Here, we observe that the PWLA \( \mathcal{P}(e^{n_t} \mid L_t,U_t,\bp^t) \), similar to \( e^{n_t} \), is a strictly increasing function in \( n_t \). At optimality, since we want to minimize \( n_t \) as much as possible, its value should satisfy the equation:
\[
n_t = \mathcal{P}^{-1}(u^t({\by,\bx}) - \xi^u_t({\by,\bx})),
\]
where \( \mathcal{P}^{-1}(z) \) is the inverse function of \( \mathcal{P}(e^{n_t} \mid L_t,U_t,\bp^t) \), which always exists, satisfying:
\[
\mathcal{P}(\exp(\mathcal{P}^{-1}(z)) \mid L_t,U_t,\bp^t) = z, \quad \forall z\in \mathbb{R}.
\]
From this observation, we can write \eqref{prob:proof-1} equivalently as:
\begin{align}
      \min_{\bx,\by}  \left\{\widehat{\cF}(\by,\bx) = \sum_{t\in [T]} \frac{u^t(\by,\bx) - \xi^u_t(\by,\bx) + \xi^t(n_t)}{h^t(\by, \bx) - \xi^h_t(\by,\bx)} \right\}\nonumber,
\end{align}
where \( n_t =  \mathcal{P}^{-1}(u^t({\by,\bx}) - \xi^u_t({\by,\bx}) ) \). In general, the approximation problem in \eqref{prob:approx-convex-exp} can be reformulated as:
\[
\min_{(\by,\bx)\in \cZ} \widehat{\cF}(\by,\bx).
\]
As a result, the approximation errors can be bounded by analyzing the gap between the original objective function \( \cF(\by,\bx) \) and its approximation \( \widehat{\cF}(\by,\bx) \). We establish this bound in the following lemma.
\begin{lemma}\label{lem:2}
    For any \( (\by,\bx) \in \cZ \), we have:
    \[
    |\widehat{\cF}(\by,\bx) - \cF(\by,\bx)| \leq  \sum_{t\in [T]} \left| \frac{2 \left( U^u_t \epsilon^h_t + U^h_t (\epsilon + \epsilon^u_t) \right)}{(L^h_t)^2 - (\epsilon^h_t)^2} \right|,
    \]
    where
    \[
    \epsilon^u_t = (\alpha L^{ht}_i + L^{gt}_i) \frac{u_i - l_i}{K}, \quad
    \epsilon^h_t =  \sum_{i\in [m]} L^{ht}_i \frac{u_i - l_i}{K},
    \]
 and \( L^u_t \) and \( U^u_t \) are the lower and upper bounds of \( u^t(\by, \bx) \), and \( L^h_t \) and \( U^h_t \) are the lower and upper bounds of \( h^t(\by, \bx) \), for all \( (\by,\bx) \in \cZ \).
\end{lemma}

With all the bounds established above, we are now ready to derive an upper bound for the approximation error yielded by any solution obtained from the approximate problem \eqref{prob:approx-convex-exp}.

\begin{theorem}\label{theorem2}
Let \( (\widehat{\by},\widehat{\bx}) \) be an optimal solution returned by solving the approximation problem \eqref{prob:approx-convex-exp}, and let \( (\by^*,\bx^*) \) be the optimal solution to the original fractional program \eqref{prob:SoR-original}. We can bound the gap between the objective values given by \( (\widehat{\by},\widehat{\bx}) \) and the optimal value as:
\[
|\cF(\widehat{\by},\widehat{\bx}) -  \cF(\by^*,\bx^*)| \leq 2\sum_{t\in [T]} \left| \frac{2 \left( U^u_t \epsilon^h_t + U^h_t (\epsilon + \epsilon^u_t) \right)}{(L^h_t)^2 - (\epsilon^h_t)^2} \right|.
\]
\end{theorem}

While the bound appears complex, we note that \( \epsilon^u_t = \mathcal{O}(1/K) \) and \( \epsilon^h_t = \mathcal{O}(1/K) \). Moreover, when \( \epsilon \) and \( 1/K \) are sufficiently small, the term \( \epsilon^h_t \) is dominated by \( L^h_t \), which further implies that the approximation error \( |\cF(\widehat{\by},\widehat{\bx}) -  \cF(\by^*,\bx^*)| \) is in \( \mathcal{O}(\epsilon + 1/K) \). This result indicates that the approximation error decreases linearly as the number of breakpoints \( K \) increases and as \( \epsilon \) approaches zero.

In practice, the approximation of the exponential function \( e^{n_t} \) can be done more efficiently than that of other univariate functions \( h^t_i(x_i) \) and \( u^t_i(x_i) \). Hence, one can choose \( \epsilon \) to be significantly small, dominated by \( 1/K \). Under this setting, the approximation error \( |\cF(\widehat{\by},\widehat{\bx}) -  \cF(\by^*,\bx^*)| \) is in \( \mathcal{O}(1/K) \).

This result formalizes the intuition that increasing the number of breakpoints $K$ will lead to improved solution accuracy. More importantly, it provides a rigorous theoretical foundation showing that the approximation error decreases at a quantifiable rate. Specifically, as $K$ increases and the exponential approximation gap $\epsilon$ decreases, our method exhibits \textit{linear convergence} toward the optimal solution of the original non-convex problem. This establishes that our approach not only scales with controllable precision but also ensures asymptotic optimality under mild regularity conditions.

\section{Applications}\label{sec:app}
We briefly discuss the applications of our discretization and approximation approach to two prominent classes of decision-making problems: the \emph{maximum capture facility location} problem and the \emph{joint assortment and price optimization} problem.

\subsection{Joint Facility Location and Cost Optimization in Maximum Capture Problem}

Let $[m]$ be the set of available locations for setting up new facilities. For each customer segment $t \in [T]$, let the utility of location $i \in [m]$ be given by
$v_{ti} = x_i \eta_{ti} + \kappa_{ti},$
where $x_i$ denotes the cost spent at location $i$, $\eta_{ti} > 0$ is a cost sensitivity parameter (reflecting how cost affects utility), and $\kappa_{ti}$ captures other utility-affecting factors such as location features. Under the logit model, the probability that a customer from segment $t$ chooses facility $i$ over competitors is given by:
\[
P^t\left(i \,\Big|\, [m] \cup \{0\}\right) = 
\frac{\exp(x_i \eta_{ti} + \kappa_{ti})}{U^t_C + \sum_{j \in [m]} \exp(x_j \eta_{tj} + \kappa_{tj})},
\]
where $0$ refers to a competitor’s facility and $U^t_C$ represents the total utility of all competing facilities.

The objective is to maximize expected captured demand (the expected number of customers attracted by the selected facilities), often referred to as the \textit{maximum capture problem} (MCP). This can be formulated as:
\begin{align}
\label{prob:maxcap}
\tag{{\sf MCP}}
\max_{(\by, \bx) \in \cZ} \left\{ f(\by, \bx) = \sum_{t \in [T]} \frac{\cQ_t \sum_{i \in [m]} y_i \exp(x_i \eta_{ti} + \kappa_{ti})}{U^t_C + \sum_{i \in [m]} y_i \exp(x_i \eta_{ti} + \kappa_{ti})} \right\},
\end{align}
where $\cQ_t$ is the proportion of customers in segment $t$. We can reformulate this into a sum-of-ratios form where decision variables only appear in the denominators:
\[
f(\by, \bx) = \sum_{t \in [T]} \cQ_t - \sum_{t \in [T]} \frac{\cQ_t U^t_C}{U^t_C + \sum_{i \in [m]} y_i \exp(x_i \eta_{ti} + \kappa_{ti})}.
\]
To apply our approximation methods, we can let $g^t_i(x_i) = 0$, $a_t = \cQ_t U^t_C$, $b_t = U^t_C$, and $h^t_i(x_i) = \exp(x_i \eta_{ti} + \kappa_{ti})$. This places the MCP within the general problem class defined by \eqref{prob:SoR-original}, allowing all of our solution methods to be applied. {\textit{ Importantly, this is the first known work to address the cost optimization aspect of MCP under random utility models with continuous cost variables.}}

\subsection{Joint Assortment and Price Optimization}

We now describe the joint assortment and price optimization problem (denoted as A\&P) under the mixed-logit model. While we reuse some notation from the previous section, the variables may take different interpretations here. Let $[m]$ denote the set of available products and $0$ denote a no-purchase option. Let $x_i$ be the price of product $i \in [m]$. For each customer segment $t \in [T]$, the utility of product $i$ is modeled as:
$v_{ti} = x_i \eta_{ti} + \kappa_{ti},$
where $\eta_{ti} < 0$ is the price sensitivity (negative, since higher prices reduce utility), and $\kappa_{ti}$ accounts for product-specific attributes. The purchase probability under the mixed-logit model is:
\[
P\left(i \,\Big|\, [m] \cup \{0\}\right) = 
\frac{\exp(x_i \eta_{ti} + \kappa_{ti})}{1 + \sum_{j \in [m]} \exp(x_j \eta_{tj} + \kappa_{tj})},
\]
where the value $1$ in the denominator represents the utility of the no-purchase option.
The goal is to jointly select an assortment and set prices to maximize expected revenue:
\begin{align}
\label{prob:asso-prive}
\tag{{\sf Assort-Price}}
\max_{\by \in \cY, \bx \in \cX} \left\{ f(\by, \bx) = \sum_{t \in [T]} \frac{\sum_{i \in [m]} y_i x_i \exp(x_i \eta_{ti} + \kappa_{ti})}{1 + \sum_{i \in [m]} y_i \exp(x_i \eta_{ti} + \kappa_{ti})} \right\}.
\end{align}
Here, $y_i$ indicates if product $i$ is offered (assortment decision) and $x_i$ is its price. To apply our methods, we let: $g^t_i(x_i) = x_i \exp(x_i \eta_{ti} + \kappa_{ti}) \text{ and } h^t_i(x_i) = \exp(x_i \eta_{ti} + \kappa_{ti})$.
Some  business  constraints can include:
\begin{itemize}
    \item Cardinality: $\sum_{i \in [m]} y_i \leq M$ (maximum number of products offered),
    \item  Price bounds: $x_i \in [l_i, u_i]$,
    \item  Budget-like constraints on bundles: $\sum_{i \in S} y_i x_i \leq W$ for a subset $S \subset [m]$.
\end{itemize}
Note that Assumption~\ref{assum:a2} holds under such constraints, enabling the use of the simplified approximation model in \eqref{prob:approx-convex-exp}. Furthermore, constraints like $\sum_{i \in S} y_i x_i \leq W$ can be easily linearized using McCormick inequalities or handled directly by Gurobi.

While joint assortment and price optimization has been widely studied, most prior work has focused on single-ratio formulations or has assumed fixed prices (see, e.g.,~\cite{wang2012capacitated,GallegoWang2014}). {\textit{In contrast, our work addresses a much more general and realistic setting by solving the joint assortment and pricing problem under the mixed-logit model with flexible linear constraints on both assortment and prices.}}

\section{Numerical Experiments}
\label{sec:experiments}
\subsection{Experimental Setup} 

We present experiments comparing our proposed approach against several baselines. To the best of our knowledge, there is currently no existing method capable of solving the general nonlinear sum-of-ratios problem \eqref{prob:SoR-original} to near-optimality with provable performance guarantees—except for general-purpose solvers such as SCIP. Therefore, we compare our approach, which combines log-transformation and PWLA, against several direct baselines discussed in Section \ref{sec:milp_socp}, as well as SCIP, a state-of-the-art solver for mixed-integer nonlinear programs.

Specifically, we consider the following approaches for comparison:
\begin{itemize}
    \item \textbf{LOG-PW }(Log-transformation + PWLA): Our proposed method, which applies a log-transformation followed by PWLA as described in Section \ref{sec:micp}. We specifically use B\&C to solve the convex program in \eqref{prob:approx-convex-exp} with valid cuts described in \eqref{ctr-OA-1}, \eqref{ctr-OA-2}, and \eqref{ctr-OA-3}, with a note that CP can provide similar performance across all instances. 
    \item \textbf{MILP }(PWLA + MILP Reformulation): The MILP-based approximation approach introduced in Section \ref{sec:milp_socp} and detailed in Supplementary Materials--Appendix C.
    \item \textbf{SOCP} (PWLA + Second-Order Cone Reformulation): The conic approximation using PWLA and second-order cone programming (SOCP), also presented in Section \ref{sec:milp_socp} and Supplementary Materials--Appendix C.
    \item \textbf{GP }(Gurobi's PWLA + Bilinear Reformulation): Gurobi's built-in PWLA capabilities combined with our bilinear reformulation, described in Section \ref{sec:milp_socp} with implementation details as in the \eqref{prob:linear-frac-Blinear} model.
    \item \textbf{SCIP}: One of the best open-source solvers for mixed-integer nonlinear programming~\citep{SCIP2024}. We use SCIP to directly solve the original nonlinear formulations \eqref{prob:maxcap} and \eqref{prob:asso-prive}. To ensure a fair comparison, we configure SCIP to use a minimum of 8 threads (as its default setting is single-threaded), with no upper limit on the number of threads for parallel computation.
\end{itemize}

To provide an overview of the formulation sizes across different approaches, Table \ref{tab:problem-static} summarizes the number of variables and constraints based on \( T \) and \( m \). It is evident that the \textbf{MILP} and \textbf{SOCP} formulations are the most complex, requiring the highest number of variables and constraints. In contrast, the formulation used in the \textbf{SCIP} solver introduces minimal additional elements, adding only one extra variable and constraint to \eqref{prob:SoR-original} to convert the objective function into a nonlinear constraint. Although the GUROBI's PWLA-based formulation appears to involve fewer variables and constraints compared to the \textbf{MILP}, \textbf{SOCP} and \textbf{LOG-PW} methods, the execution of the built-in function \texttt{addGenConstrExp()} introduces several new binary variables and constraints. This results in an overall formulation that can become significantly larger than that  of the \textbf{MILP} or \textbf{SOCP} approaches.

\begin{table}[htb]
\centering
\resizebox{\textwidth}{!}{%
\begin{tabular}{l|ccc|ccc}
\multicolumn{1}{l|}{} &
  \multicolumn{3}{c|}{MCP} &
  \multicolumn{3}{c}{A\&P} \\ \cline{2-7} 
\multicolumn{1}{l|}{} &
  \multicolumn{2}{c|}{\#Variables} &
  \multirow{2}{*}{\#Constraints} &
  \multicolumn{2}{c|}{\#Variables} &
  \multirow{2}{*}{\#Constraints} \\ \cline{2-3} \cline{5-6}
\multicolumn{1}{l|}{} &
  \multicolumn{1}{c|}{Binary} &
  \multicolumn{1}{c|}{Continuous} &
   &
  \multicolumn{1}{c|}{Binary} &
  \multicolumn{1}{c|}{Continuous} &
   \\ \hline
\multicolumn{1}{r|}{MILP} &
  \multicolumn{1}{c|}{$m+mK$} &
  \multicolumn{1}{c|}{$T+m+Tm+TmK$} &
  $T+2m+mK+4Tm+4TmK+2$ &
  \multicolumn{1}{c|}{$m+mK$} &
  \multicolumn{1}{c|}{$T+2m+mK+Tm+TmK$} &
  $T+2m+mK+4Tm+4TmK+2$ \\ \hline
\multicolumn{1}{r|}{SOCP} &
  \multicolumn{1}{c|}{$m+mK$} &
  \multicolumn{1}{c|}{$2T+m+Tm+TmK$} &
  $3T+2m+mK+Tm+TmK+2$ &
  \multicolumn{1}{c|}{$m+mK$} &
  \multicolumn{1}{c|}{$2T+m+Tm+TmK$} &
  $3T+2m+mK +5Tm+5TmK+2$ \\ \hline
\multicolumn{1}{r|}{GP} &
  \multicolumn{1}{c|}{$m$} &
  \multicolumn{1}{c|}{$3T+m+2Tm$} &
  $3T+2Tm+2$ &
  \multicolumn{1}{c|}{$m$} &
  \multicolumn{1}{c|}{$3T+m+3Tm$} &
  $3T+m+6Tm+2$ \\ \hline
\multicolumn{1}{r|}{SCIP} &
  \multicolumn{1}{c|}{$m$} &
  \multicolumn{1}{c|}{$m+1$} &
  3 &
  \multicolumn{1}{c|}{$m$} &
  \multicolumn{1}{c|}{$m+1$} &
  3 \\ \hline
  \multicolumn{1}{r|}{LOG-PW} &
  \multicolumn{1}{c|}{$m + mK + \sum_{t=1}^T H_t$} &
  \multicolumn{1}{c|}{$T + m + \sum_{t=1}^T H_t$} &
  $mK + 5m + 3\sum_{t=1}^T H_t + \mathcal{C}$ &
  \multicolumn{1}{c|}{$m + mK + \sum_{t=1}^T H_t$} &
  \multicolumn{1}{c|}{$4T + m + \sum_{t=1}^T H_t$} &
  $2T + mK + 4m + 3\sum_{t=1}^T H_t + \mathcal{C}$ \\ \hline
\end{tabular}%
}
\caption{Problem formulation sizes used in \textbf{MILP}, \textbf{SOCP}, \textbf{GP}, \textbf{SCIP} and \textbf{LOG-PW} approaches; $\mathcal{C}$ is number of lazy constraints added to the model when applying GUROBI's \texttt{callback()} functions; $H_t$ is number of breakpoints used to linearize $e^{n_t}$.}
\label{tab:problem-static}
\end{table}


The experiments are conducted on a PC with processors Intel(R) Core(TM) i7-9700 CPU @ 3.00GHz,
RAM of 16 gigabytes, and operating system Window 11. The code is in C++
and links to GUROBI 11.0.3 (under default settings) to solve the \textbf{MILP}, \textbf{SOCP}, \textbf{LOG-PW} and SCIP version 9.1.1 for the \textbf{SCIP} model. We set the CPU time limit for
each instance as 3600 seconds, i.e., we stop the algorithms/solver if they exceed the time budget and report the best solutions found. We provide numerical experiments based on two applications: A\&P and MCP.  We use the same number of PWLA segments for both our methods (\textbf{LOG-PW, MILP, SOCP}) and the \textbf{GP} approach.

\subsection{{Solution Quality as $K$ and $\epsilon$ Change}}

\paragraph{Choice of $K$.} In this experiment, we analyze the performance of our PWLA approach as a function of $K$, aiming to determine the best choice for achieving high-quality solutions in practice. Since the size of the reformulations in \eqref{prob:approx-convex-exp} scales proportionally with $K$, selecting an appropriate value is essential to balancing computational cost and solution accuracy. The preliminary experiment show that for $K \geq 25$, the objective gaps become negligible, suggesting that $K=25$ is a reasonable choice to achieve near-optimal performance across most instances (see Figures \ref{fig:K_T_2} and \ref{fig:K_T_5} in Appendix \ref{appd:K_tau} for details). Based on these findings, we fix $K = 25$ in all subsequent experiments.

\paragraph{Choice of $\epsilon$.}
In the approximation model \eqref{prob:approx-convex-exp}, the parameter $\epsilon$ governs the granularity of the PWLA for the exponential function $e^{n_t}$. This parameter is critical in controlling the trade-off between approximation accuracy and computational efficiency. Following our preliminary experiment (shown in Figures \ref{fig:e_T_2} and \ref{fig:e_T_5} in Appendix \ref{appd:K_tau}), we select $\epsilon = 10^{-3}$ for all subsequent experiments, as it provides a robust balance between approximation accuracy and computational runtime. This choice ensures that our \textbf{LOG-PW} approach remains both efficient and reliable across all problem instances considered.

\subsection{Comparison Results}\label{sec:result}
We now present a comprehensive comparison between our proposed method (\textbf{LOG-PW}) and all the aforementioned baselines across various instances of both the MCP and the A\&P problems. The following performance metrics are reported:
\begin{itemize}
    \item \textbf{Optimality Count:} The number of instances for which each method solves the problem to optimality. For PWLA-based methods, this refers to solving the approximated model optimally. For SCIP, it refers to solving the original nonlinear problem to proven optimality within the time limit.
    \item \textbf{Best Objective Count:} For each method, we take the solution it returns and evaluate its objective value under the original formulation \eqref{prob:SoR-original}. We then count how many times each method achieves the best objective value compared to all other baselines.
    \item \textbf{Average Runtime:} The average computation time (in seconds) across instances solved to near-optimality. Since we set a time limit of 3600 seconds, if the average runtime exceeds this budget, we indicate it with ``-'' in the results table.
\end{itemize}

\begin{table}[h]\footnotesize
\centering
\resizebox{\textwidth}{!}{%
\begin{tabular}{llll|rrrrr|rrrrr|rrrrr}
 &
   &
   &
   &
  \multicolumn{5}{c|}{\#Instances  solved optimally} &
  \multicolumn{5}{c|}{\#Instances with best objectives} &
  \multicolumn{5}{c}{Average time (s)} \\ \hline
\multicolumn{1}{l|}{$T$} &
  $m$ &
  $C$ &
  $M$ &
  MILP &
  SOCP &
  GP &
  SCIP &
  LOG-PW &
  MILP &
  SOCP &
  GP &
  SCIP &
  LOG-PW &
  MILP &
  SOCP &
  GP &
  SCIP &
  LOG-PW \\ \hline
\multicolumn{1}{l|}{\multirow{8}{*}{5}} &
  \multirow{4}{*}{50} &
  \multirow{2}{*}{20} &
  16 &
  \textbf{3} &
  \textbf{3} &
  \textbf{3} &
  0 &
  \textbf{3} &
  \textbf{3} &
  \textbf{3} &
  \textbf{3} &
  2 &
  \textbf{3} &
  67.24 &
  0.21 &
  0.21 &
  - &
  \textbf{0.08} \\
\multicolumn{1}{l|}{} &
   &
   &
  25 &
  \textbf{3} &
  \textbf{3} &
  \textbf{3} &
  0 &
  \textbf{3} &
  \textbf{3} &
  \textbf{3} &
  \textbf{3} &
  0 &
  \textbf{3} &
  69.46 &
  0.18 &
  0.86 &
  - &
  \textbf{0.08} \\
\multicolumn{1}{l|}{} &
   &
  \multirow{2}{*}{30} &
  16 &
  \textbf{3} &
  \textbf{3} &
  \textbf{3} &
  0 &
  \textbf{3} &
  \textbf{3} &
  \textbf{3} &
  \textbf{3} &
  2 &
  \textbf{3} &
  62.67 &
  0.18 &
  0.23 &
  - &
  \textbf{0.08} \\
\multicolumn{1}{l|}{} &
   &
   &
  25 &
  \textbf{3} &
  \textbf{3} &
  \textbf{3} &
  0 &
  \textbf{3} &
  \textbf{3} &
  \textbf{3} &
  \textbf{3} &
  0 &
  \textbf{3} &
  42.25 &
  0.18 &
  1.04 &
  - &
  \textbf{0.08} \\ \cline{2-19} 
\multicolumn{1}{l|}{} &
  \multirow{4}{*}{100} &
  \multirow{2}{*}{40} &
  33 &
  1 &
  \textbf{3} &
  \textbf{3} &
  0 &
  \textbf{3} &
  \textbf{3} &
  \textbf{3} &
  \textbf{3} &
  0 &
  \textbf{3} &
  3162 &
  0.27 &
  1.29 &
  - &
  \textbf{0.12} \\
\multicolumn{1}{l|}{} &
   &
   &
  50 &
  0 &
  \textbf{3} &
  \textbf{3} &
  0 &
  \textbf{3} &
  \textbf{3} &
  \textbf{3} &
  \textbf{3} &
  1 &
  \textbf{3} &
  - &
  0.28 &
  1.53 &
  - &
  \textbf{0.12} \\
\multicolumn{1}{l|}{} &
   &
  \multirow{2}{*}{60} &
  33 &
  0 &
  \textbf{3} &
  \textbf{3} &
  0 &
  \textbf{3} &
  \textbf{3} &
  \textbf{3} &
  \textbf{3} &
  0 &
  \textbf{3} &
  - &
  0.29 &
  1.83 &
  - &
  \textbf{0.11} \\
\multicolumn{1}{l|}{} &
   &
   &
  50 &
  1 &
  \textbf{3} &
  \textbf{3} &
  0 &
  \textbf{3} &
  \textbf{3} &
  \textbf{3} &
  \textbf{3} &
  0 &
  \textbf{3} &
  2495.01 &
  0.3 &
  3.06 &
  - &
  \textbf{0.11} \\ \hline
\multicolumn{1}{l|}{\multirow{8}{*}{10}} &
  \multirow{4}{*}{50} &
  \multirow{2}{*}{20} &
  16 &
  0 &
  \textbf{3} &
  \textbf{3} &
  0 &
  \textbf{3} &
  \textbf{3} &
  \textbf{3} &
  \textbf{3} &
  0 &
  \textbf{3} &
  - &
  0.28 &
  1.46 &
  - &
  0.11 \\
\multicolumn{1}{l|}{} &
   &
   &
  25 &
  0 &
  \textbf{3} &
  \textbf{3} &
  0 &
  \textbf{3} &
  \textbf{3} &
  \textbf{3} &
  \textbf{3} &
  0 &
  \textbf{3} &
  - &
  0.27 &
  1.92 &
  - &
  \textbf{0.11} \\
\multicolumn{1}{l|}{} &
   &
  \multirow{2}{*}{30} &
  16 &
  0 &
  \textbf{3} &
  \textbf{3} &
  0 &
  \textbf{3} &
  \textbf{3} &
  \textbf{3} &
  \textbf{3} &
  0 &
  \textbf{3} &
  - &
  0.26 &
  0.40 &
  - &
  \textbf{0.10} \\
\multicolumn{1}{l|}{} &
   &
   &
  25 &
  0 &
  \textbf{3} &
  \textbf{3} &
  0 &
  \textbf{3} &
  \textbf{3} &
  \textbf{3} &
  \textbf{3} &
  0 &
  \textbf{3} &
  - &
  0.27 &
  1.38 &
  - &
  \textbf{0.11} \\ \cline{2-19} 
\multicolumn{1}{l|}{} &
  \multirow{4}{*}{100} &
  \multirow{2}{*}{40} &
  33 &
  0 &
  \textbf{3} &
  \textbf{3} &
  0 &
  \textbf{3} &
  \textbf{3} &
  \textbf{3} &
  \textbf{3} &
  0 &
  \textbf{3} &
  - &
  0.49 &
  4.23 &
  - &
  \textbf{0.20} \\
\multicolumn{1}{l|}{} &
   &
   &
  50 &
  0 &
  \textbf{3} &
  \textbf{3} &
  0 &
  \textbf{3} &
  \textbf{3} &
  \textbf{3} &
  \textbf{3} &
  0 &
  \textbf{3} &
  - &
  0.50 &
  10.63 &
  - &
  \textbf{0.16} \\
\multicolumn{1}{l|}{} &
   &
  \multirow{2}{*}{60} &
  33 &
  0 &
  \textbf{3} &
  \textbf{3} &
  0 &
  \textbf{3} &
  \textbf{3} &
  \textbf{3} &
  \textbf{3} &
  1 &
  \textbf{3} &
  - &
  0.53 &
  5.90 &
  - &
  \textbf{0.19} \\
\multicolumn{1}{l|}{} &
   &
   &
  50 &
  0 &
  \textbf{3} &
  \textbf{3} &
  0 &
  \textbf{3} &
  \textbf{3} &
  \textbf{3} &
  \textbf{3} &
  0 &
  \textbf{3} &
  - &
  0.51 &
  14.69 &
  - &
  \textbf{0.28} \\ \hline
\multicolumn{1}{l|}{\multirow{8}{*}{100}} &
  \multirow{4}{*}{50} &
  \multirow{2}{*}{20} &
  16 &
  0 &
  \textbf{3} &
  0 &
  0 &
  \textbf{3} &
  \textbf{3} &
  \textbf{3} &
  \textbf{3} &
  0 &
  \textbf{3} &
  - &
  7.23 &
  - &
  - &
  \textbf{0.76} \\
\multicolumn{1}{l|}{} &
   &
   &
  25 &
  0 &
  \textbf{3} &
  0 &
  0 &
  \textbf{3} &
  \textbf{3} &
  \textbf{3} &
  \textbf{3} &
  0 &
  \textbf{3} &
  - &
  5.35 &
  - &
  - &
  \textbf{0.65} \\
\multicolumn{1}{l|}{} &
   &
  \multirow{2}{*}{30} &
  16 &
  0 &
  \textbf{3} &
  0 &
  0 &
  \textbf{3} &
  \textbf{3} &
  \textbf{3} &
  \textbf{3} &
  0 &
  \textbf{3} &
  - &
  7.12 &
  - &
  - &
  \textbf{0.63} \\
\multicolumn{1}{l|}{} &
   &
   &
  25 &
  0 &
  \textbf{3} &
  0 &
  0 &
  \textbf{3} &
  \textbf{3} &
  \textbf{3} &
  \textbf{3} &
  0 &
  \textbf{3} &
  - &
  5.5 &
  - &
  - &
  \textbf{0.64} \\ \cline{2-19} 
\multicolumn{1}{l|}{} &
  \multirow{4}{*}{100} &
  \multirow{2}{*}{40} &
  33 &
  0 &
  \textbf{3} &
  0 &
  0 &
  \textbf{3} &
  \textbf{3} &
  \textbf{3} &
  \textbf{3} &
  0 &
  \textbf{3} &
  - &
  11.83 &
  - &
  - &
  \textbf{1.71} \\
\multicolumn{1}{l|}{} &
   &
   &
  50 &
  0 &
  \textbf{3} &
  0 &
  0 &
  \textbf{3} &
  \textbf{3} &
  \textbf{3} &
  \textbf{3} &
  0 &
  \textbf{3} &
  - &
  12.96 &
  - &
  - &
  \textbf{1.45} \\
\multicolumn{1}{l|}{} &
   &
  \multirow{2}{*}{60} &
  33 &
  0 &
  \textbf{3} &
  0 &
  0 &
  \textbf{3} &
  \textbf{3} &
  \textbf{3} &
  \textbf{3} &
  0 &
  \textbf{3} &
  - &
  11.97 &
  - &
  - &
  \textbf{1.51} \\
\multicolumn{1}{l|}{} &
   &
   &
  50 &
  0 &
  \textbf{3} &
  0 &
  0 &
  \textbf{3} &
  \textbf{3} &
  \textbf{3} &
  \textbf{3} &
  0 &
  \textbf{3} &
  - &
  11.98 &
  - &
  - &
  \textbf{1.47} \\ \hline
\multicolumn{4}{r|}{Summary:} &
  14 &
  \textbf{72} &
  48 &
  0 &
  \textbf{72} &
  \textbf{72} &
  \textbf{72} &
  \textbf{72} &
  6 &
  \textbf{72} &
   &
   &
   &
   &
   \\ \hline
\end{tabular}%
}
\caption{Comparison results for MCP instances of large $T$; instances  are grouped by $(T,m,C,M)$.}
\label{table3}
\end{table}

\begin{table}[!h]\footnotesize
\centering
\resizebox{0.7\textwidth}{!}{%
\begin{tabular}{rrr|rrr|rrr|rrr}
\multicolumn{1}{l}{} &
  \multicolumn{1}{l}{} &
  \multicolumn{1}{l|}{} &
  \multicolumn{3}{c|}{\#Solved optimally} &
  \multicolumn{3}{c|}{\#Best objective} &
  \multicolumn{3}{c}{Average time (s)} \\ \hline
\multicolumn{1}{r}{$m$} &
  $C$ &
  $M$ &
  SOCP &
  GP &
  LOG-PW &
  SOCP &
  GP &
  LOG-PW &
  SOCP &
  GP &
  LOG-PW \\ \hline
\multicolumn{1}{r}{\multirow{4}{*}{50}} &
  \multirow{2}{*}{20} &
  16 &
  \textbf{3} & \textbf{3} & \textbf{3} & \textbf{3} & \textbf{3} & \textbf{3} & 0.28 & 1.46 & \textbf{0.11} \\
\multicolumn{1}{r}{} &
   &
  25 &
  \textbf{3} & \textbf{3} & \textbf{3} & \textbf{3} & \textbf{3} & \textbf{3} & 0.27 & 1.92 & \textbf{0.11} \\
\multicolumn{1}{r}{} &
  \multirow{2}{*}{30} &
  16 &
  \textbf{3} & \textbf{3} & \textbf{3} & \textbf{3} & \textbf{3} & \textbf{3} & 0.26 & 0.40 & \textbf{0.10} \\
\multicolumn{1}{r}{} &
   &
  25 &
  \textbf{3} & \textbf{3} & \textbf{3} & \textbf{3} & \textbf{3} & \textbf{3} & 0.27 & 1.38 & \textbf{0.11} \\ \hline
\multicolumn{1}{r}{\multirow{4}{*}{100}} &
  \multirow{2}{*}{40} &
  33 &
  \textbf{3} & \textbf{3} & \textbf{3} & \textbf{3} & \textbf{3} & \textbf{3} & 0.49 & 4.23 & \textbf{0.20} \\
\multicolumn{1}{r}{} &
   &
  50 &
  \textbf{3} & \textbf{3} & \textbf{3} & \textbf{3} & \textbf{3} & \textbf{3} & 0.50 & 10.63 & \textbf{0.16} \\
\multicolumn{1}{r}{} &
  \multirow{2}{*}{60} &
  33 &
  \textbf{3} & \textbf{3} & \textbf{3} & \textbf{3} & \textbf{3} & \textbf{3} & 0.53 & 5.90 & \textbf{0.19} \\
\multicolumn{1}{r}{} &
   &
  50 &
  \textbf{3} & \textbf{3} & \textbf{3} & \textbf{3} & \textbf{3} & \textbf{3} & 0.51 & 14.69 & \textbf{0.28} \\ \hline
\multicolumn{1}{r}{\multirow{4}{*}{200}} &
  \multirow{2}{*}{80} &
  66 &
  \textbf{3} & \textbf{3} & \textbf{3} & \textbf{3} & \textbf{3} & \textbf{3} & 1.08 & 10.74 & \textbf{0.30} \\
\multicolumn{1}{r}{} &
   &
  100 &
  \textbf{3} & \textbf{3} & \textbf{3} & \textbf{3} & \textbf{3} & \textbf{3} & 1.07 & 63.74 & \textbf{0.30} \\
\multicolumn{1}{r}{} &
  \multirow{2}{*}{120} &
  66 &
  \textbf{3} & \textbf{3} & \textbf{3} & \textbf{3} & \textbf{3} & \textbf{3} & 1.07 & 27.03 & \textbf{0.40} \\
\multicolumn{1}{r}{} &
   &
  100 &
  \textbf{3} & \textbf{3} & \textbf{3} & \textbf{3} & \textbf{3} & \textbf{3} & 1.24 & 14.22 & \textbf{0.34} \\ \hline
\multicolumn{1}{r}{\multirow{4}{*}{500}} &
  \multirow{2}{*}{200} &
  166 &
  \textbf{3} & 2 & \textbf{3} & \textbf{3} & \textbf{3} & \textbf{3} & 8.35 & 2673.31 & \textbf{0.80} \\
\multicolumn{1}{r}{} &
   &
  250 &
  \textbf{3} & 2 & \textbf{3} & \textbf{3} & \textbf{3} & \textbf{3} & 6.30 & 2882.27 & \textbf{0.80} \\
\multicolumn{1}{r}{} &
  \multirow{2}{*}{300} &
  166 &
  \textbf{3} & 1 & \textbf{3} & \textbf{3} & \textbf{3} & \textbf{3} & 6.56 & 617.88 & \textbf{0.82} \\
\multicolumn{1}{r}{} &
   &
  250 &
  \textbf{3} & 1 & \textbf{3} & \textbf{3} & \textbf{3} & \textbf{3} & 8.78 & 2209.59 & \textbf{0.77} \\ \hline
\multicolumn{1}{r}{\multirow{4}{*}{1000}} &
  \multirow{2}{*}{400} &
  333 &
  \textbf{3} & 0 & \textbf{3} & \textbf{3} & \textbf{3} & \textbf{3} & 34.66 & - & \textbf{2.01} \\
\multicolumn{1}{r}{} &
   &
  500 &
  \textbf{3} & \textbf{3} & \textbf{3} & \textbf{3} & \textbf{3} & \textbf{3} & 37.68 & 147.02 & \textbf{1.92} \\
\multicolumn{1}{r}{} &
  \multirow{2}{*}{600} &
  333 &
  \textbf{3} & 0 & \textbf{3} & \textbf{3} & \textbf{3} & \textbf{3} & 29.16 & - & \textbf{1.72} \\
\multicolumn{1}{r}{} &
   &
  500 &
  \textbf{3} & \textbf{3} & \textbf{3} & \textbf{3} & \textbf{3} & \textbf{3} & 31.03 & 156.14 & \textbf{1.78} \\ \hline
\multicolumn{3}{r|}{Summary:} &
  \textbf{60} & 48 & \textbf{60} & \textbf{60} & \textbf{60} & \textbf{60} &  &  &  \\ \hline
\end{tabular}%
}
\caption{Comparison results for MCP instances of large $m$ and $T=10$; instances are grouped by $(m,C,M)$.}
\label{table3_extend}
\end{table}

\paragraph{MCP Instances.}

We begin with the joint facility location and cost optimization problem, i.e., the MCP. To make the problem more realistic, we incorporate two constraints: (i) a cardinality constraint on the number of selected locations, represented as $\sum_{i \in [m]} y_i \leq M$, and (ii) an upper bound on the total cost spent on opening new facilities, given by $\sum_{i \in [m]} y_i x_i \leq C$. We conduct numerical comparisons across five different approaches: \textbf{LOG-PW}, \textbf{MILP}, \textbf{SOCP}, \textbf{GP}, and \textbf{SCIP} (which directly solves the original nonlinear formulation). For each setting defined by the tuple $(T, m, C, M)$, we generate 3 independent instances and solve them using all five methods. 

The comparison results are shown in Tables \ref{table3} and \ref{table3_extend}, where the best results are shown in bold. In Table \ref{table3}, we present the results for instances with a large value of $T$ (up to 100), while maintaining a small or medium number of locations $m$. In contrast, Table \ref{table3_extend} reports the results for instances with a large number of locations ($m$ varies from 200 to 1000), while keeping $T$ at a small value of 10. Note that the \textbf{MILP} and  \textbf{SOCP} solvers cannot solve any instance with $T=10$, therefore, the results of these methods are not included in Table \ref{table3_extend}.

It is not surprising to see that \textbf{MILP} and \textbf{GP} are outperformed by the other methods, in terms of solution quality. The two best approaches are \textbf{LOG-PW} and \textbf{SOCP}, respectively. These methods solve all instances to optimal, and the \textbf{LOG-PW} provides the shorter runtime than \textbf{SOCP}. Interestingly, \textbf{SOCP} clearly outperforms \textbf{MILP} in terms of both solution quality and computing time --- \textbf{SOCP} is able to return best objective values for all the instances and the maximum running time is just about 37.68 seconds, while \textbf{MILP} cannot return the best objectives for several large-sized instances and always exceeds the time budget of 3600 seconds.


\begin{table}[!h]\footnotesize
\centering
\resizebox{\textwidth}{!}{%
\begin{tabular}{rrrr|rrrrr|rrrrr|rrrrr}
 &
   &
   &
   &
  \multicolumn{5}{c|}{\#Solved optimally} &
  \multicolumn{5}{c|}{\#Best objective} &
  \multicolumn{5}{c}{Average runtime (s)} \\ \hline
\multicolumn{1}{r|}{$T$} &
  $m$ &
  $C$ &
  $M$ &
  MILP &
  SOCP &
  GP &
  SCIP &
  LOG-PW &
  MILP &
  SOCP &
  GP &
  SCIP &
  LOG-PW &
  MILP &
  SOCP &
  GP &
  SCIP &
  LOG-PW\\ \hline
\multicolumn{1}{r|}{\multirow{16}{*}{2}} &
  \multirow{4}{*}{10} &
  \multirow{2}{*}{4} &
  3 &
  \textbf{3} &
  \textbf{3} &
  \textbf{3} &
  \textbf{3} &
  \textbf{3} &
  1 &
  1 &
  \textbf{3} &
  \textbf{3} &
  1 &
  1.48 &
  1.63 &
  \textbf{0.22} &
  3.92 &
  2.42 \\
\multicolumn{1}{r|}{} &
   &
   &
  5 &
  \textbf{3} &
  \textbf{3} &
  \textbf{3} &
  2 &
  \textbf{3} &
  1 &
  1 &
  \textbf{3} &
  \textbf{3} &
  1 &
  2.31 &
  2.09 &
  \textbf{0.55} &
  242.24 &
  4.73 \\
\multicolumn{1}{r|}{} &
   &
  \multirow{2}{*}{6} &
  3 &
  \textbf{3} &
  \textbf{3} &
  \textbf{3} &
  \textbf{3} &
  \textbf{3} &
  0 &
  0 &
  \textbf{3} &
  \textbf{3} &
  0 &
  1.79 &
  2.32 &
  \textbf{0.20} &
  10.16 &
  2.27 \\
\multicolumn{1}{r|}{} &
   &
   &
  5 &
  \textbf{3} &
  \textbf{3} &
  \textbf{3} &
  \textbf{3} &
  \textbf{3} &
  0 &
  0 &
  \textbf{3} &
  \textbf{3} &
  0 &
  3.56 &
  3.87 &
  \textbf{0.95} &
  932.15 &
  6.13 \\ \cline{2-19} 
\multicolumn{1}{r|}{} &
  \multirow{4}{*}{20} &
  \multirow{2}{*}{8} &
  6 &
  \textbf{3} &
  \textbf{3} &
  \textbf{3} &
  1 &
  \textbf{3} &
  2 &
  2 &
  \textbf{3} &
  \textbf{3} &
  2 &
  28.96 &
  32.83 &
  \textbf{1.30} &
  3388.39 &
  7.60 \\
\multicolumn{1}{r|}{} &
   &
   &
  10 &
  \textbf{3} &
  \textbf{3} &
  \textbf{3} &
  0 &
  \textbf{3} &
  2 &
  2 &
  \textbf{3} &
  1 &
  2 &
  135.42 &
  174.53 &
  \textbf{1.23} &
  - &
  7.01 \\
\multicolumn{1}{r|}{} &
   &
  \multirow{2}{*}{12} &
  6 &
  \textbf{3} &
  \textbf{3} &
  \textbf{3} &
  0 &
  \textbf{3} &
  2 &
  2 &
  2 &
  \textbf{3} &
  2 &
  39.17 &
  38.50 &
  \textbf{1.55} &
  - &
  8.94 \\
\multicolumn{1}{r|}{} &
   &
   &
  10 &
  \textbf{3} &
  \textbf{3} &
  \textbf{3} &
  0 &
  \textbf{3} &
  \textbf{3} &
  \textbf{3} &
  2 &
  0 &
  \textbf{3} &
  515.15 &
  226.25 &
  \textbf{2.86} &
  - &
  7.34 \\ \cline{2-19} 
\multicolumn{1}{r|}{} &
  \multirow{4}{*}{50} &
  \multirow{2}{*}{20} &
  16 &
  0 &
  0 &
  \textbf{3} &
  0 &
  \textbf{3} &
  \textbf{3} &
  2 &
  \textbf{3} &
  0 &
  \textbf{3} &
  - &
  - &
  \textbf{7.41} &
  - &
  25.84 \\
\multicolumn{1}{r|}{} &
   &
   &
  25 &
  0 &
  0 &
  \textbf{3} &
  0 &
  \textbf{3} &
  2 &
  1 &
  \textbf{3} &
  0 &
  \textbf{3} &
  - &
  - &
  \textbf{7.47} &
  - &
  16.15 \\
\multicolumn{1}{r|}{} &
   &
  \multirow{2}{*}{30} &
  16 &
  0 &
  0 &
  \textbf{3} &
  0 &
  \textbf{3} &
  \textbf{3} &
  \textbf{3} &
  1 &
  0 &
  \textbf{3} &
  - &
  - &
  36.28 &
  - &
  \textbf{16.06} \\
\multicolumn{1}{r|}{} &
   &
   &
  25 &
  0 &
  0 &
  \textbf{3} &
  0 &
  \textbf{3} &
  \textbf{3} &
  \textbf{3} &
  1 &
  0 &
  \textbf{3} &
  - &
  - &
  \textbf{16.65} &
  - &
  40.32 \\ \cline{2-19} 
\multicolumn{1}{r|}{} &
  \multirow{4}{*}{100} &
  \multirow{2}{*}{40} &
  33 &
  0 &
  0 &
  2 &
  0 &
  \textbf{3} &
  1 &
  0 &
  \textbf{3} &
  0 &
  \textbf{3} &
  - &
  - &
  2441.79 &
  - &
  \textbf{111.13} \\
\multicolumn{1}{r|}{} &
   &
   &
  50 &
  0 &
  0 &
  1 &
  0 &
  \textbf{3} &
  2 &
  0 &
  \textbf{3} &
  0 &
  \textbf{3} &
  - &
  - &
  2442.00 &
  - &
  \textbf{70.19} \\
\multicolumn{1}{r|}{} &
   &
  \multirow{2}{*}{60} &
  33 &
  0 &
  0 &
  0 &
  0 &
  \textbf{3} &
  2 &
  \textbf{3} &
  0 &
  0 &
  \textbf{3} &
  - &
  - &
  - &
  - &
  \textbf{69.83} \\
\multicolumn{1}{r|}{} &
   &
   &
  50 &
  0 &
  0 &
  0 &
  0 &
  \textbf{3} &
  2 &
  1 &
  0 &
  0 &
  \textbf{3} &
  - &
  - &
  - &
  - &
  \textbf{64.13} \\ \hline
\multicolumn{1}{r|}{\multirow{16}{*}{5}} &
  \multirow{4}{*}{10} &
  \multirow{2}{*}{4} &
  3 &
  \textbf{3} &
  \textbf{3} &
  \textbf{3} &
  \textbf{3} &
  \textbf{3} &
  1 &
  1 &
  2 &
  \textbf{3} &
  1 &
  4.32 &
  4.66 &
  \textbf{0.29} &
  118.09 &
  9.00 \\
\multicolumn{1}{r|}{} &
   &
   &
  5 &
  \textbf{3} &
  \textbf{3} &
  \textbf{3} &
  2 &
  \textbf{3} &
  1 &
  1 &
  2 &
  \textbf{3} &
  1 &
  12.91 &
  11.60 &
  \textbf{0.41} &
  1407.28 &
  10.59 \\
\multicolumn{1}{r|}{} &
   &
  \multirow{2}{*}{6} &
  3 &
  \textbf{3} &
  \textbf{3} &
  \textbf{3} &
  \textbf{3} &
  \textbf{3} &
  0 &
  0 &
  \textbf{3} &
  \textbf{3} &
  0 &
  4.77 &
  4.97 &
  \textbf{0.34} &
  213.39 &
  11.23 \\
\multicolumn{1}{r|}{} &
   &
   &
  5 &
  \textbf{3} &
  \textbf{3} &
  \textbf{3} &
  1 &
  \textbf{3} &
  0 &
  0 &
  2 &
  \textbf{3} &
  0 &
  10.29 &
  9.30 &
  \textbf{0.62} &
  3507.58 &
  16.38 \\ \cline{2-19} 
\multicolumn{1}{r|}{} &
  \multirow{4}{*}{20} &
  \multirow{2}{*}{8} &
  6 &
  \textbf{3} &
  \textbf{3} &
  \textbf{3} &
  0 &
  \textbf{3} &
  2 &
  2 &
  1 &
  \textbf{3} &
  2 &
  272.58 &
  284.3 &
  \textbf{3.91} &
  - &
  32.98 \\
\multicolumn{1}{r|}{} &
   &
   &
  10 &
  2 &
  2 &
  \textbf{3} &
  0 &
  \textbf{3} &
  \textbf{2} &
  1 &
  \textbf{2} &
  0 &
  \textbf{2} &
  2893.25 &
  2879.31 &
  \textbf{7.49} &
  - &
  110.84 \\
\multicolumn{1}{r|}{} &
   &
  \multirow{2}{*}{12} &
  6 &
  \textbf{3} &
  \textbf{3} &
  \textbf{3} &
  0 &
  \textbf{3} &
  2 &
  2 &
  2 &
  \textbf{3} &
  2 &
  142.63 &
  147.38 &
  \textbf{5.62} &
  - &
  27.48 \\
\multicolumn{1}{r|}{} &
   &
   &
  10 &
  2 &
  2 &
  \textbf{3} &
  0 &
  \textbf{3} &
  \textbf{3} &
  2 &
  2 &
  0 &
  \textbf{3} &
  2324.7 &
  2780.4 &
  \textbf{13.79} &
  - &
  38.05 \\ \cline{2-19} 
\multicolumn{1}{r|}{} &
  \multirow{4}{*}{50} &
  \multirow{2}{*}{20} &
  16 &
  0 &
  0 &
  2 &
  0 &
  \textbf{3} &
  \textbf{3} &
  0 &
  \textbf{3} &
  0 &
  \textbf{3} &
  - &
  - &
  1357.45 &
  - &
  \textbf{521.46} \\
\multicolumn{1}{r|}{} &
   &
   &
  25 &
  0 &
  0 &
  \textbf{3} &
  0 &
  \textbf{3} &
  1 &
  0 &
  \textbf{3} &
  0 &
  \textbf{3} &
  - &
  - &
  1358.93 &
  - &
  \textbf{237.84} \\
\multicolumn{1}{r|}{} &
   &
  \multirow{2}{*}{30} &
  16 &
  0 &
  0 &
  0 &
  0 &
  \textbf{3} &
  \textbf{3} &
  0 &
  0 &
  0 &
  \textbf{3} &
  - &
  - &
  - &
  - &
  \textbf{85.53} \\
\multicolumn{1}{r|}{} &
   &
   &
  25 &
  0 &
  0 &
  0 &
  0 &
  \textbf{3} &
  2 &
  0 &
  0 &
  0 &
  \textbf{3} &
  - &
  - &
  - &
  - &
  \textbf{209.55} \\ \cline{2-19} 
\multicolumn{1}{r|}{} &
  \multirow{4}{*}{100} &
  \multirow{2}{*}{40} &
  33 &
  0 &
  0 &
  0 &
  0 &
  \textbf{3} &
  1 &
  0 &
  1 &
  1 &
  \textbf{3} &
  - &
  - &
  - &
  - &
  \textbf{493.81} \\
\multicolumn{1}{r|}{} &
   &
   &
  50 &
  0 &
  0 &
  0 &
  0 &
  \textbf{3} &
  0 &
  0 &
  1 &
  0 &
  \textbf{3} &
  - &
  - &
  - &
  - &
  \textbf{1924.06} \\
\multicolumn{1}{r|}{} &
   &
  \multirow{2}{*}{60} &
  33 &
  0 &
  0 &
  0 &
  0 &
  \textbf{3} &
  0 &
  0 &
  0 &
  0 &
  \textbf{3} &
  - &
  - &
  - &
  - &
  \textbf{354.02} \\
\multicolumn{1}{r|}{} &
   &
   &
  50 &
  0 &
  0 &
  0 &
  0 &
  \textbf{3} &
  0 &
  0 &
  0 &
  0 &
  \textbf{3} &
  - &
  - &
  - &
  - &
  \textbf{932.89} \\ \hline
\multicolumn{4}{r|}{Summary:} &
  46 &
  46 &
  68 &
  21 &
  \textbf{96} &
  51 &
  32 &
  60 &
  38 &
  \textbf{70} &
   &
   &
   &
   &
   \\ \hline
\end{tabular}%
}
\caption{Comparison results for A\&P instances with $T=2$ and $5$; instances grouped by $(T,m,C,M)$.}
\label{table1_extend}
\end{table}

\begin{table}[!h]\footnotesize
\centering
\resizebox{\textwidth}{!}{%
\begin{tabular}{rrr|rrrrrr|rrrrrr}
\multicolumn{1}{l}{} &
  \multicolumn{1}{l}{} &
  \multicolumn{1}{l|}{} &
  \multicolumn{6}{c|}{$T = 10$} &
  \multicolumn{6}{c}{$T = 20$} \\ \cline{4-15} 
\multicolumn{1}{l}{} &
  \multicolumn{1}{l}{} &
  \multicolumn{1}{l|}{} &
  \multicolumn{2}{c|}{\#Solved optimally} &
  \multicolumn{2}{c|}{\#Best objective} &
  \multicolumn{2}{c|}{Average runtime (s)} &
  \multicolumn{2}{c|}{\#Solved optimally} &
  \multicolumn{2}{c|}{\#Best objective} &
  \multicolumn{2}{c}{Average runtime (s)} \\ \hline
\multicolumn{1}{r}{$m$} &
  $C$ &
  $M$ &
  GP &
  \multicolumn{1}{r|}{LOG-PW} &
  GP &
  \multicolumn{1}{r|}{LOG-PW} &
  GP &
  LOG-PW &
  GP &
  \multicolumn{1}{r|}{LOG-PW} &
  GP &
  \multicolumn{1}{r|}{LOG-PW} &
  GP &
  LOG-PW \\ \hline
\multicolumn{1}{r}{\multirow{4}{*}{10}} &
  \multirow{2}{*}{4} &
  3 &
  \textbf{3} &
  \multicolumn{1}{r|}{\textbf{3}} &
  \textbf{3} &
  \multicolumn{1}{r|}{0} &
  10.75 &
  \textbf{4.64} &
  \textbf{3} &
  \multicolumn{1}{r|}{\textbf{3}} &
  \textbf{3} &
  \multicolumn{1}{r|}{1} &
  68.74 &
  \textbf{10.91} \\
\multicolumn{1}{r}{} &
   &
  5 &
  \textbf{3} &
  \multicolumn{1}{r|}{\textbf{3}} &
  \textbf{3} &
  \multicolumn{1}{r|}{0} &
  40.41 &
  \textbf{14.82} &
  \textbf{3} &
  \multicolumn{1}{r|}{\textbf{3}} &
  \textbf{3} &
  \multicolumn{1}{r|}{1} &
  368.6 &
  \textbf{52.18} \\
\multicolumn{1}{r}{} &
  \multirow{2}{*}{6} &
  3 &
  \textbf{3} &
  \multicolumn{1}{r|}{\textbf{3}} &
  \textbf{3} &
  \multicolumn{1}{r|}{2} &
  \textbf{0.83} &
  7.50 &
  \textbf{3} &
  \multicolumn{1}{r|}{\textbf{3}} &
  \textbf{3} &
  \multicolumn{1}{r|}{1} &
  17.97 &
  \textbf{3.71} \\
\multicolumn{1}{r}{} &
   &
  5 &
  \textbf{3} &
  \multicolumn{1}{r|}{\textbf{3}} &
  \textbf{3} &
  \multicolumn{1}{r|}{1} &
  44.34 &
  \textbf{5.13} &
  \textbf{3} &
  \multicolumn{1}{r|}{\textbf{3}} &
  \textbf{3} &
  \multicolumn{1}{r|}{1} &
  836.38 &
  \textbf{15.19} \\ \hline
\multicolumn{1}{r}{\multirow{4}{*}{20}} &
  \multirow{2}{*}{8} &
  6 &
  0 &
  \multicolumn{1}{r|}{\textbf{3}} &
  \textbf{3} &
  \multicolumn{1}{r|}{\textbf{3}} &
  - &
  \textbf{23.41} &
  0 &
  \multicolumn{1}{r|}{\textbf{3}} &
  2 &
  \multicolumn{1}{r|}{\textbf{3}} &
  - &
  \textbf{342.29} \\
\multicolumn{1}{r}{} &
   &
  10 &
  0 &
  \multicolumn{1}{r|}{\textbf{3}} &
  1 &
  \multicolumn{1}{r|}{\textbf{3}} &
  - &
  \textbf{21.72} &
  0 &
  \multicolumn{1}{r|}{\textbf{3}} &
  0 &
  \multicolumn{1}{r|}{\textbf{3}} &
  - &
  \textbf{1012.63} \\
\multicolumn{1}{r}{} &
  \multirow{2}{*}{12} &
  6 &
  \textbf{3} &
  \multicolumn{1}{r|}{\textbf{3}} &
  \textbf{3} &
  \multicolumn{1}{r|}{\textbf{3}} &
  93.50 &
  \textbf{9.81} &
  0 &
  \multicolumn{1}{r|}{\textbf{3}} &
  2 &
  \multicolumn{1}{r|}{\textbf{3}} &
  - &
  \textbf{28.05} \\
\multicolumn{1}{r}{} &
   &
  10 &
  1 &
  \multicolumn{1}{r|}{\textbf{3}} &
  \textbf{2} &
  \multicolumn{1}{r|}{\textbf{2}} &
  2402.45 &
  \textbf{59.70} &
  0 &
  \multicolumn{1}{r|}{\textbf{3}} &
  0 &
  \multicolumn{1}{r|}{\textbf{3}} &
  - &
  \textbf{2139.88} \\ \hline
\multicolumn{1}{r}{\multirow{4}{*}{50}} &
  \multirow{2}{*}{20} &
  16 &
  0 &
  \multicolumn{1}{r|}{\textbf{3}} &
  0 &
  \multicolumn{1}{r|}{\textbf{3}} &
  - &
  \textbf{535.41} &
  0 &
  \multicolumn{1}{r|}{\textbf{3}} &
  0 &
  \multicolumn{1}{r|}{\textbf{3}} &
  - &
  \textbf{35.29} \\
\multicolumn{1}{r}{} &
   &
  25 &
  0 &
  \multicolumn{1}{r|}{\textbf{3}} &
  0 &
  \multicolumn{1}{r|}{\textbf{3}} &
  - &
  \textbf{193.36} &
  0 &
  \multicolumn{1}{r|}{\textbf{3}} &
  0 &
  \multicolumn{1}{r|}{\textbf{3}} &
  - &
  \textbf{1624.7} \\
\multicolumn{1}{r}{} &
  \multirow{2}{*}{30} &
  16 &
  0 &
  \multicolumn{1}{r|}{\textbf{3}} &
  0 &
  \multicolumn{1}{r|}{\textbf{3}} &
  - &
  \textbf{43.93} &
  0 &
  \multicolumn{1}{r|}{\textbf{3}} &
  0 &
  \multicolumn{1}{r|}{\textbf{3}} &
  - &
  \textbf{94.24} \\
\multicolumn{1}{r}{} &
   &
  25 &
  0 &
  \multicolumn{1}{r|}{\textbf{3}} &
  0 &
  \multicolumn{1}{r|}{\textbf{3}} &
  - &
  \textbf{165.4} &
  0 &
  \multicolumn{1}{r|}{\textbf{3}} &
  0 &
  \multicolumn{1}{r|}{\textbf{3}} &
  - &
  \textbf{310.28} \\ \hline
\multicolumn{1}{r}{\multirow{4}{*}{100}} &
  \multirow{2}{*}{40} &
  33 &
  0 &
  \multicolumn{1}{r|}{\textbf{3}} &
  0 &
  \multicolumn{1}{r|}{\textbf{3}} &
  - &
  \textbf{969.99} &
  0 &
  \multicolumn{1}{r|}{\textbf{3}} &
  0 &
  \multicolumn{1}{r|}{\textbf{3}} &
  - &
  \textbf{577.49} \\
\multicolumn{1}{r}{} &
   &
  50 &
  0 &
  \multicolumn{1}{r|}{\textbf{3}} &
  0 &
  \multicolumn{1}{r|}{\textbf{3}} &
  - &
  \textbf{1193.56} &
  0 &
  \multicolumn{1}{r|}{1} &
  0 &
  \multicolumn{1}{r|}{\textbf{3}} &
  - &
  \textbf{2968.3} \\
\multicolumn{1}{r}{} &
  \multirow{2}{*}{60} &
  33 &
  0 &
  \multicolumn{1}{r|}{\textbf{3}} &
  0 &
  \multicolumn{1}{r|}{\textbf{3}} &
  - &
  \textbf{221.33} &
  0 &
  \multicolumn{1}{r|}{\textbf{3}} &
  0 &
  \multicolumn{1}{r|}{\textbf{3}} &
  - &
  \textbf{176.46} \\
\multicolumn{1}{r}{} &
   &
  50 &
  0 &
  \multicolumn{1}{r|}{\textbf{3}} &
  0 &
  \multicolumn{1}{r|}{\textbf{3}} &
  - &
  \textbf{434.37} &
  0 &
  \multicolumn{1}{r|}{2} &
  0 &
  \multicolumn{1}{r|}{\textbf{3}} &
  - &
  \textbf{2032.19} \\ \hline
\multicolumn{3}{r|}{Summary:} &
  16 &
  \multicolumn{1}{r|}{\textbf{48}} &
  21 &
  \multicolumn{1}{r|}{\textbf{38}} &
   &
   &
  12 &
  \multicolumn{1}{r|}{\textbf{45}} &
  16 &
  \multicolumn{1}{r|}{\textbf{40}} &
   &
   \\ \hline
\end{tabular}%
}
\caption{Comparison results for A\&P instances with $T=10$ and $20$; instances grouped by $(T,m,C,M)$.}
\label{table2_extend}
\end{table}

\paragraph{A\&P Instances.} Let us now shift our attention to the A\&P problem. These A\&P instances pose a significantly greater challenge in terms of solving, especially when compared to the MCP ones. This heightened complexity primarily arises from the non-convex nature of the fractional program even when the continuous variables are fixed \citep{rusmevichientong2014assortment}. We, therefore, adopt a small value of $T$ small, specifically setting it to 2 and 5, while varying the number of products $m$ up to 100. Similar to the MCP instances, we introduce two constraints, i.e., a cardinality constraint on the size of the selected assortment $\sum_{i\in [m]}y_i\leq M$, and an upper bound constraint on a weighted sum of the prices $\sum_{i\in [m]} \alpha_i y_ix_i \leq C$, where $\alpha_i$ take random values in $[0.5, 1]$.  The second constraint can be described as one that mandates the total price of a given set of offered products to remain below a specified upper limit. For each group of $(m,C,M)$, we randomly generate 3 instances and report the number of instances that are solved to optimality.

In Table \ref{table1_extend}, it is clear that \textbf{LOG-PW} outperforms \textbf{GP} in both terms of solution quality and runtime. The \textbf{MILP} and \textbf{SOCP} perform worse than GUROBI's PWLA (i.e. \textbf{GP}) in returning good solutions. In cases of \textbf{MILP} and \textbf{SOCP}, large number of McCormick inequalities limits their performances when $m \geq 50$. The \textbf{GP} is the second fastest approach to solving instances with $m \leq 50$, however, it cannot handle instances with large $m$ because of the increase in complexity within the built-in approximation process, as we mentioned at the beginning of this section. 

Table \ref{table2_extend} presents the results of the two best approaches -- \textbf{LOG-PW} and \textbf{GP}, on larger datasets with $T \in \{10, 20\}$. We can see that \textbf{LOG-PW} is able to provide optimal solutions for 93 out of 96 instances, while \textbf{GP} only solves to optimal 28 instances. The \textbf{LOG-PW} also finds 78 best objective values, compared to 37 ones of the \textbf{GP}. This once again confirms the superiority of \textbf{LOG-PW} over other methods as the number of fractions in the objective function increases. The average value of the objective function deviation and the average runtime of the A\&P instances are detailed in Supplementary Materials--Appendix B.2.



\section{Conclusion}\label{sec:concl}
We studied a class of non-convex binary-continuous sum-of-ratios programs that arise in several important decision-making applications, including assortment and price optimization, as well as maximum capture facility location. To address the computational challenges posed by the nonlinearity and fractional nature of these problems, we proposed a novel and innovative solution framework based on a combination of log-transformation and PWLA. This transformation enables the reformulation of the original nonlinear fractional program into a mixed-integer convex program, where standard optimization techniques such as CP or B\&C can be applied efficiently using gradient-based valid cuts. We also established theoretical performance guarantees for the solutions obtained from the approximated model and provided practical guidance for selecting the discretization parameter $K$ to ensure near-optimality. Through extensive numerical experiments on both assortment and price optimization, and facility location and cost optimization problems, we demonstrated the effectiveness of our proposed approximation method. In particular, the \textbf{LOG-PW} approach showed superior performance compared to several baselines, including: Gurobi’s built-in PWLA, PWLA combined with \textbf{MILP} and \textbf{SOCP} reformulations, and the general-purpose mixed-integer nonlinear solver \textbf{SCIP}.

Future research directions include extending our methodology to accommodate broader classes of discrete choice models, such as the nested and cross-nested logit models~\citep{Trai03}, or the network-based Generalized Extreme Value (GEV) models~\citep{DalyBier06, MaiFreFosBas15_DynMEV}.

\bibliographystyle{apacite}
\bibliography{refs}

\clearpage

\appendix

\section*{APPENDIX}
Appendix \ref{appd:proofs} provides technical proofs that were omitted from the main paper. Appendix \ref{appd:results} provides additional experiments.

\section{Proofs}\label{appd:proofs}

\subsection{Proof of Lemma \ref{lem:1}}
    Let \( \widehat{u}^t_i(x_i) \) and \( \widehat{h}^t_i(x_i) \) be the PWLAs of \( u^t_i(x_i) \) and \( h^t_i(x_i) \), respectively. Recall that \( g^t_i(x_i) \) and \( h^t_i(x_i) \) are Lipschitz continuous with constants \( L^{gt}_i \) and \( L^{ht}_i \), respectively (Assumption A1). Moreover, since 
    \[
    u^t_i(x_i) = \alpha h^t_i(x_i) - g^t_i(x_i), \quad \forall i\in [m],
    \]
    the function \( u^t_i(x_i) \) is also Lipschitz continuous with constant \( \alpha L^{ht}_i + L^{gt}_i \). This Lipschitz continuity implies:
    \[
        |h^t_i(x_i) - \widehat{h}^t_i(x_i)| \leq  L^{ht}_i \frac{u_i - l_i}{K}, \quad \text{and} \quad
        |u^t_i(x_i) - \widehat{u}^t_i(x_i)| \leq  (\alpha L^{ht}_i + L^{gt}_i) \frac{u_i - l_i}{K}. 
    \]
    Consequently, we obtain:
    \begin{align}
        |u^t(\by,\bx) - \widehat{u}^t(\by,\bx)| &\leq \sum_{i\in [m]} y_i | u^t_i(x_i) - \widehat{u}^t_i(x_i)| \leq \sum_{i\in [m]} (\alpha L^{ht}_i + L^{gt}_i) \frac{u_i - l_i}{K}, \nonumber\\
        |h^t(\by,\bx) - \widehat{h}^t(\by,\bx)| &\leq \sum_{i\in [m]} y_i | h^t_i(x_i) - \widehat{h}^t_i(x_i)| \leq \sum_{i\in [m]} L^{ht}_i \frac{u_i - l_i}{K}, \nonumber
    \end{align}
    which establishes the desired bounds.
\qed
\subsection{Proof of Lemma \ref{lem:2}}
    From Lemma \ref{lem:1}, we have \( |\xi^u_t(\by,\bx)| \leq \epsilon^u_t \) and \( |\xi^h_t(\by,\bx)| \leq \epsilon^h_t \). To bound the gap between \( \widehat{\cF}(\by, \bx) \) and \( \cF(\by, \bx) \), we derive the following inequalities:
    \begin{align}
        \sum_{t\in [T]} \frac{u^t(\by,\bx) - (\epsilon + \epsilon^u_t)}{h^t(\by, \bx) + \epsilon^h_t} &\leq \cF(\by,\bx) \leq \sum_{t\in [T]} \frac{u^t(\by,\bx) + (\epsilon + \epsilon^u_t)}{h^t(\by, \bx) - \epsilon^h_t}, \nonumber \\
        \sum_{t\in [T]} \frac{u^t(\by,\bx) - (\epsilon + \epsilon^u_t)}{h^t(\by, \bx) + \epsilon^h_t} &\leq \widehat{\cF}(\by,\bx) \leq \sum_{t\in [T]} \frac{u^t(\by,\bx) + (\epsilon + \epsilon^u_t)}{h^t(\by, \bx) - \epsilon^h_t}. \nonumber
    \end{align}
    Thus, we obtain:
    \begin{align}
        |\cF(\by,\bx) - \widehat{\cF}(\by,\bx)| &\leq \sum_{t\in [T]} \left|\frac{u^t(\by,\bx) + (\epsilon + \epsilon^u_t)}{h^t(\by, \bx) - \epsilon^h_t}  - \frac{u^t(\by,\bx) - (\epsilon + \epsilon^u_t)}{h^t(\by, \bx) + \epsilon^h_t} \right| \nonumber\\
        &=   \sum_{t\in [T]} \left| \frac{2 \left( u^t(\by,\bx) \epsilon^h_t + h^t(\by,\bx) (\epsilon + \epsilon^u_t) \right)}{(h^t(\by,\bx))^2 - (\epsilon^h_t)^2} \right| 
        \leq \sum_{t\in [T]} \left| \frac{2 \left( U^u_t \epsilon^h_t + U^h_t (\epsilon + \epsilon^u_t) \right)}{(L^h_t)^2 - (\epsilon^h_t)^2} \right|,\nonumber
    \end{align}
    where \( L^u_t, U^u_t, L^h_t, U^h_t \) are the lower and upper bounds of \( u^t(\by, \bx) \) and \( h^t(\by, \bx) \) for all \( (\by,\bx) \in \cZ \). This completes the proof.
\qed

\subsection{Proof of Theorem \ref{theorem1}}
To show that \( \mathcal{P}(e^x \mid U, L, \bp) \) is monotonically increasing, we need to prove that for any \( x_1, x_2 \in [L, U] \) with \( x_1 < x_2 \), we have \( \mathcal{P}(e^{x_1} \mid U, L, \bp) < \mathcal{P}(e^{x_2} \mid U, L, \bp) \).

We consider two cases:
\begin{itemize}
    \item \textbf{Case 1:} If \( x_1 \) and \( x_2 \) belong to different sub-intervals, there exists a breakpoint \( p_h \) such that \( x_1 < p_h \leq x_2 \). By the definition of PWLA,
    \[
    \mathcal{P}(e^{x_2} \mid U, L, \bp) \geq e^{p_h} > \mathcal{P}(e^{x_1} \mid U, L, \bp).
    \]
    \item \textbf{Case 2:} If \( x_1, x_2 \) belong to the same sub-interval, assume \( x_1, x_2 \in [p_h, p_{h+1}] \) for some index \( h \in [H] \). The function \( \mathcal{P}(e^x \mid U, L, \bp) \) is a linear function connecting the two points \( (p_h, e^{p_h}) \) and \( (p_{h+1}, e^{p_{h+1}}) \), given by:
    \begin{align*}
        \mathcal{P}(e^{x_1} \mid U, L, \bp) &= e^{p_h} + (x_1 - p_h) \frac{e^{p_{h+1}} - e^{p_h}}{p_{h+1} - p_h}, \\
        \mathcal{P}(e^{x_2} \mid U, L, \bp) &= e^{p_h} + (x_2 - p_h) \frac{e^{p_{h+1}} - e^{p_h}}{p_{h+1} - p_h}.
    \end{align*}
    Since \( x_1 < x_2 \) and \( \frac{e^{p_{h+1}} - e^{p_h}}{p_{h+1} - p_h} > 0 \), it follows that \( \mathcal{P}(e^{x_1} \mid U, L, \bp) < \mathcal{P}(e^{x_2} \mid U, L, \bp) \).
\end{itemize}
This validates the monotonicity of \( \mathcal{P}(e^x \mid U, L, \bp) \). The existence of a well-defined inverse function \( \mathcal{P}^{-1}(z, \bp) \) follows directly from this monotonicity property.

For Case 2, from the way we select the breakpoints, we seek each next point \( p_{h+1} \) as far as possible while ensuring that the gap does not exceed \( \epsilon \). If \( p_{h+1} \) is not the upper bound \( U \), then the maximum gap should be equal to \( \epsilon \), i.e.,
\begin{equation}\label{eq:cP-eq-1}
    \max_{x \in [p_h, p_{h+1}]} \left\{ \phi(x) = e^{p_h} + (x - p_h) \frac{e^{p_{h+1}} - e^{p_h}}{p_{h+1} - p_h} - e^x \right\} = \epsilon.
\end{equation}
Moreover, we can upper-bound the gap function \( \phi(x) \) as follows:
\begin{align*}
    e^{p_h} + (x - p_h) \frac{e^{p_{h+1}} - e^{p_h}}{p_{h+1} - p_h} - e^x 
    &\leq e^{p_h} + (p_{h+1} - p_h) \frac{e^{p_{h+1}} - e^{p_h}}{p_{h+1} - p_h} - e^{p_h} = e^{p_{h+1}} - e^{p_h}.
\end{align*}
From the Mean Value Theorem, there exists \( c \in [p_h, p_{h+1}] \) such that:
\begin{equation}\label{eq:cP-eq-2}
    e^{p_{h+1}} - e^{p_h} = (p_{h+1} - p_h) e^c \leq e^U (p_{h+1} - p_h).
\end{equation}
Combining \eqref{eq:cP-eq-1} and \eqref{eq:cP-eq-2}, we obtain:
\[
e^U (p_{h+1} - p_h) \geq \epsilon, \quad \text{or equivalently,} \quad p_{h+1} - p_h \geq \frac{\epsilon}{e^U}.
\]
Summing over all intervals, we get:
\[
U - L \geq p_H - p_1 = \sum_{h=1}^{H-1} (p_{h+1} - p_h) \geq \frac{\epsilon (H-1)}{e^U}.
\]
Rearranging, we obtain the bound:
\[
H \leq \frac{e^U (U - L)}{\epsilon} + 1.
\]
This completes the proof.
\qed

\subsection{Proof of Theorem \ref{theorem2}}
    For notational simplicity, let us define:
    \[
    C = \sum_{t\in [T]} \left| \frac{2 \left( U^u_t \epsilon^h_t + U^h_t (\epsilon + \epsilon^u_t) \right)}{(L^h_t)^2 - (\epsilon^h_t)^2} \right|.
    \]
    Since both the original fractional program and the approximate problem \eqref{prob:approx-convex-exp} share the same feasible set, we have:
    \begin{align*}
        \cF(\widehat{\by},\widehat{\bx}) &\geq \cF(\by^*,\bx^*)
        \stackrel{(a)}{\geq} \widehat{\cF}(\by^*,\bx^*) - C
        \stackrel{(b)}{\geq} \widehat{\cF}(\widehat{\by},\widehat{\bx}) - C
        \stackrel{(c)}{\geq} \cF(\widehat{\by},\widehat{\bx}) - 2C.
    \end{align*}
    Here, inequalities \( (a) \) and \( (c) \) follow from Lemma \ref{lem:2}, while \( (b) \) holds because \( (\widehat{\by},\widehat{\bx}) \) is optimal for the approximate problem with objective \( \widehat{\cF}(\by,\bx) \).     This directly implies:
    \[
    \cF(\widehat{\by},\widehat{\bx}) - \cF(\by^*,\bx^*) \leq 2C,
    \]
    as desired.
\qed

\section{Additional Experiment Results}\label{appd:results}

\subsection{On the Choice of $K$ and $\epsilon$}\label{appd:K_tau}
To evaluate the solution quality for different values of $K$ and $\epsilon$, we set the value of $\tau$ equal to $10^{-4}$ then solve all assortment and pricing (A\&P) instances with $T = 2$ and $T = 5$ (see Section \ref{sec:result} for details on instance generation). For each instance, we vary $K \in \{5, 10, 15, 20, 25, 30, 35, 40, 45, 50\}$ with a fixed approximation threshold $\epsilon = 10^{-6}$, which controls the approximation gap for the exponential function $e^{n_t}$ as discussed in Section \ref{sec:e_nt}. We solve the corresponding approximate problem for each $K$ and obtain the solution denoted by $(\by^K, \bx^K)$. To estimate the value of $\epsilon$, the \textbf{LOG-PW} is run with $\epsilon \in \{10^{-6}, 10^{-5},...,10^{-2}\}$ and $K$ found by the above process to archive all solutions $(\by^\epsilon, \bx^\epsilon)$.

To assess solution quality, we compute the percentage gap between the objective value $f(\by^K, \bx^K)$, $f(\by^\epsilon, \bx^\epsilon)$ and the objective value returned by \textbf{SCIP}. It is worth noting that the solutions from \textbf{SCIP} are not necessarily optimal, but are used as a common baseline for comparison across different values of $K$ and $\epsilon$. A positive gap indicates that our method (\textbf{LOG-PW}) yields a better solution, while a negative gap indicates that SCIP performs better.

\begin{figure}[!h]
    \centering
    \begin{subfigure}[b]{0.48\textwidth}
        \centering
        \includegraphics[width=\textwidth]{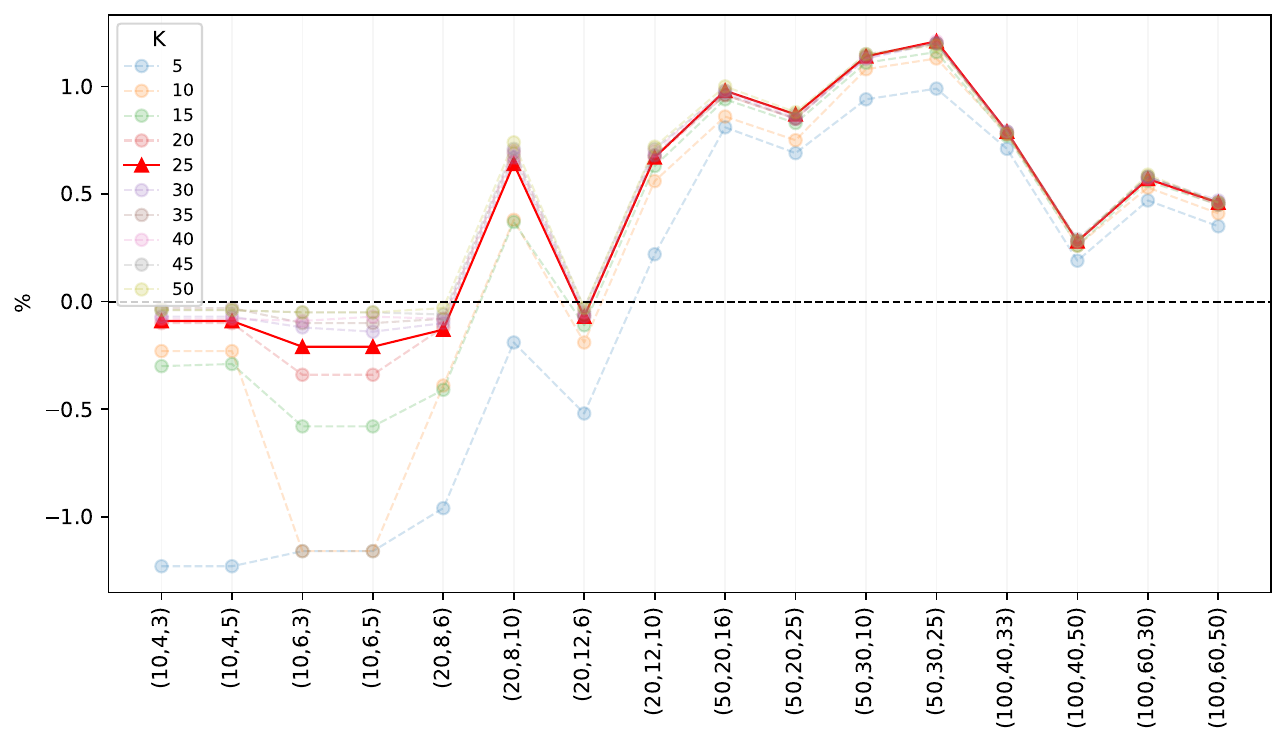} 
        \caption{Average gap of objective value}
        \label{fig:obj_T_2}
    \end{subfigure}
    \begin{subfigure}[b]{0.48\textwidth}
        \centering
        \includegraphics[width=\textwidth]{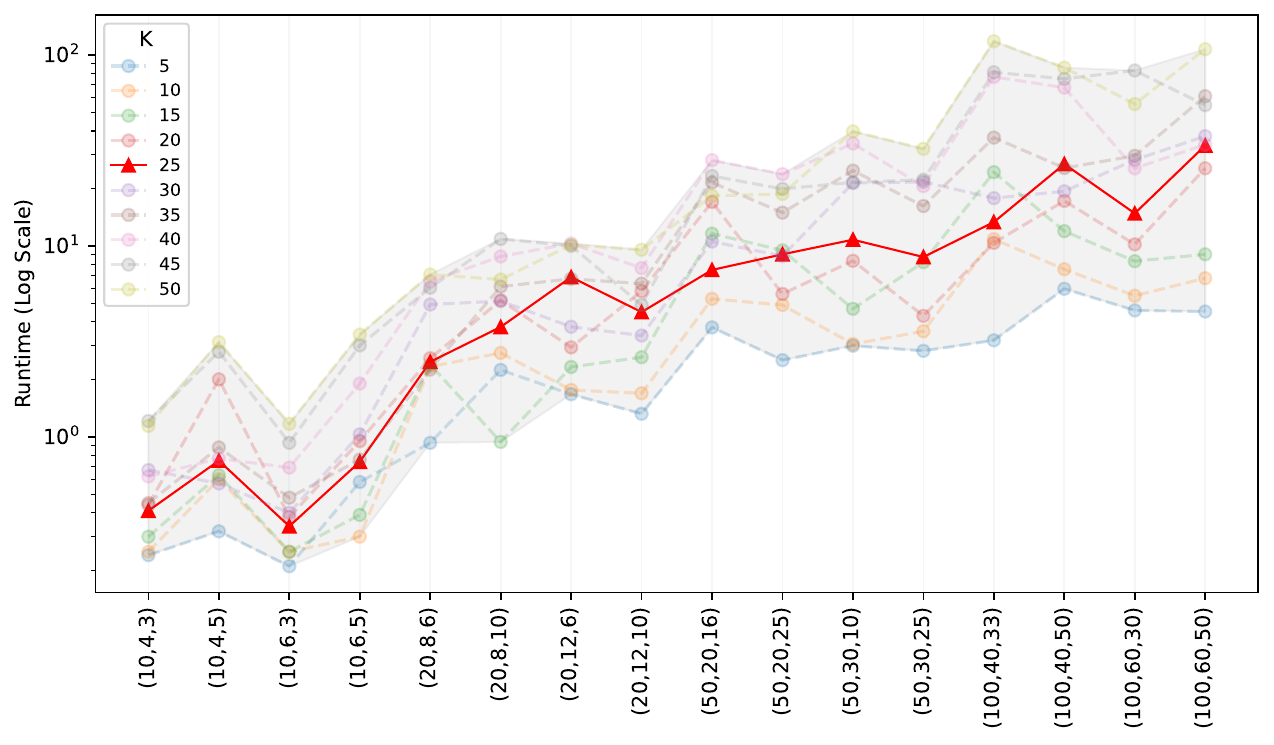} 
        \caption{Average runtime}
        \label{fig:time_T_2}
    \end{subfigure}
    
    \caption{Impact of $K$ to \textbf{LOG-PW} on A\&P instances with $T=2$}
    \label{fig:K_T_2}
\end{figure}

\begin{figure}[!h]
    \centering
    \begin{subfigure}[b]{0.48\textwidth}
        \centering
        \includegraphics[width=\textwidth]{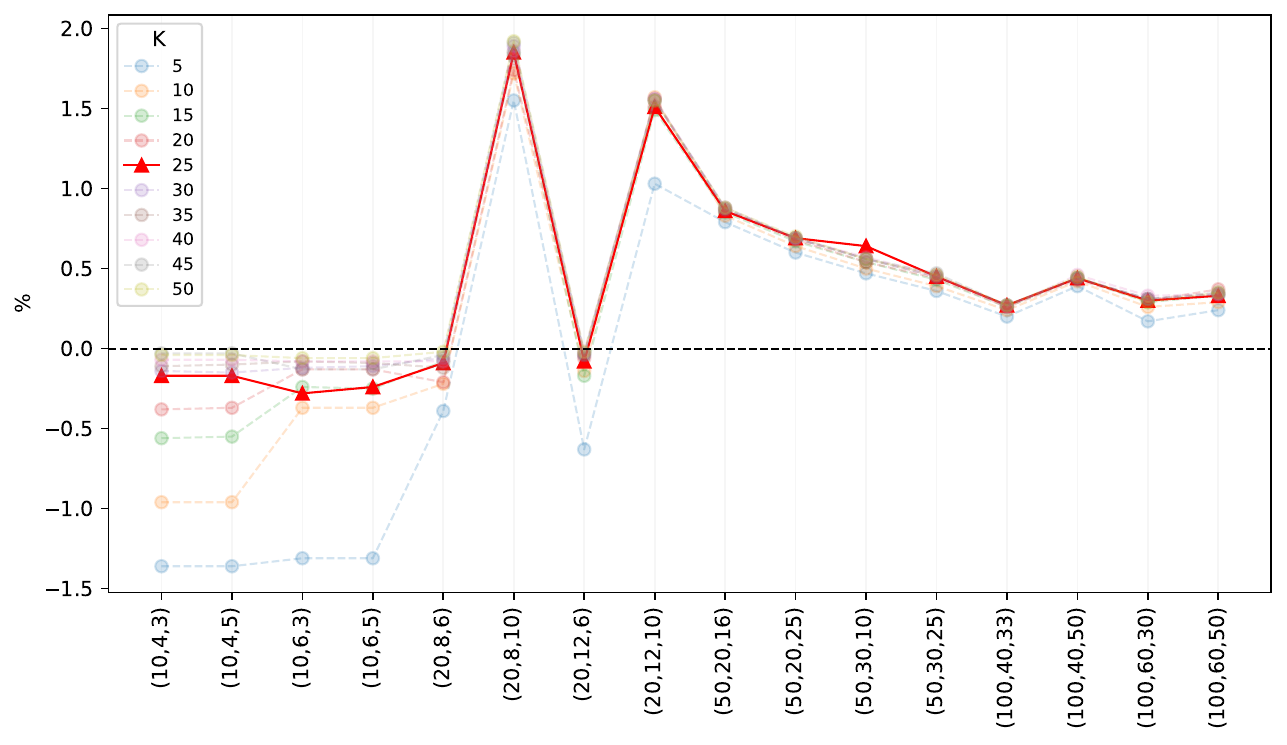} 
        \caption{Average gap of objective value}
        \label{fig:obj_T_5}
    \end{subfigure}
    \begin{subfigure}[b]{0.48\textwidth}
        \centering
        \includegraphics[width=\textwidth]{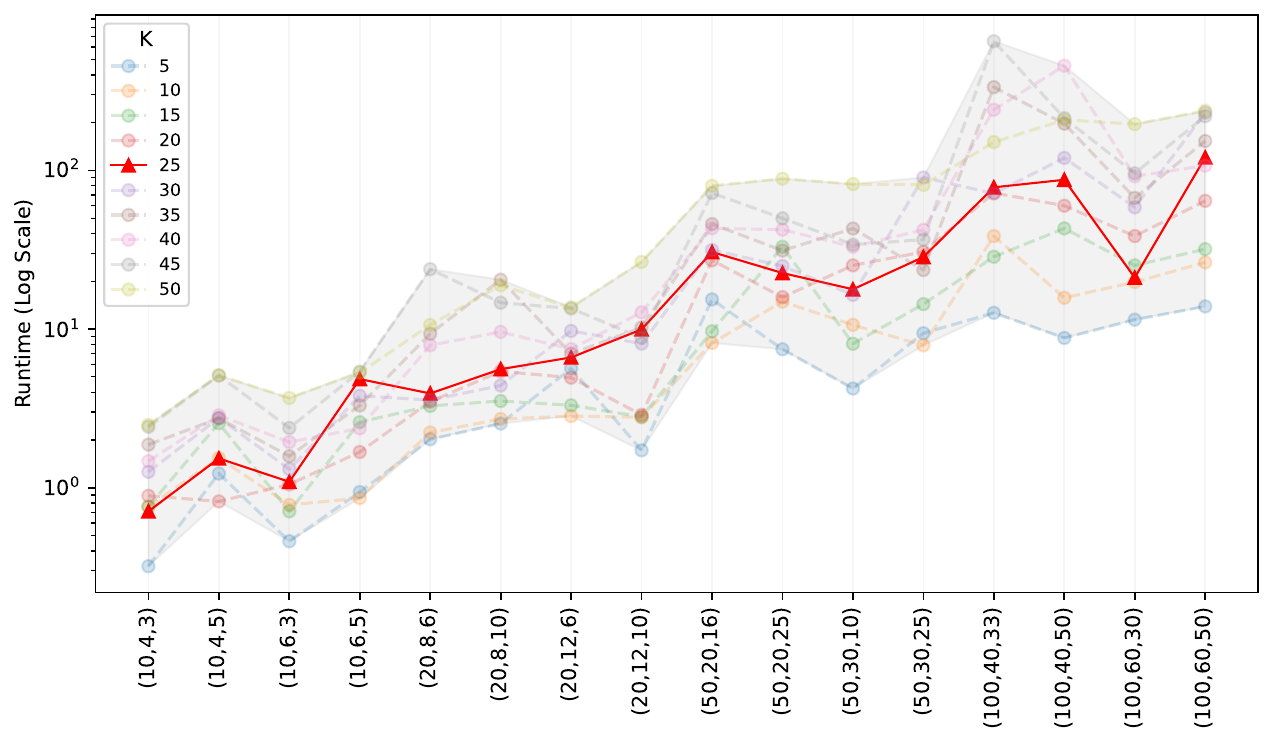} 
        \caption{Average runtime}
        \label{fig:time_T_5}
    \end{subfigure}
    
    \caption{Impact of $K$ to \textbf{LOG-PW} on A\&P instances with $T=5$}
    \label{fig:K_T_5}
\end{figure}

\begin{figure}[!h]
    \centering
    \begin{subfigure}[b]{0.48\textwidth}
        \centering
        \includegraphics[width=\textwidth]{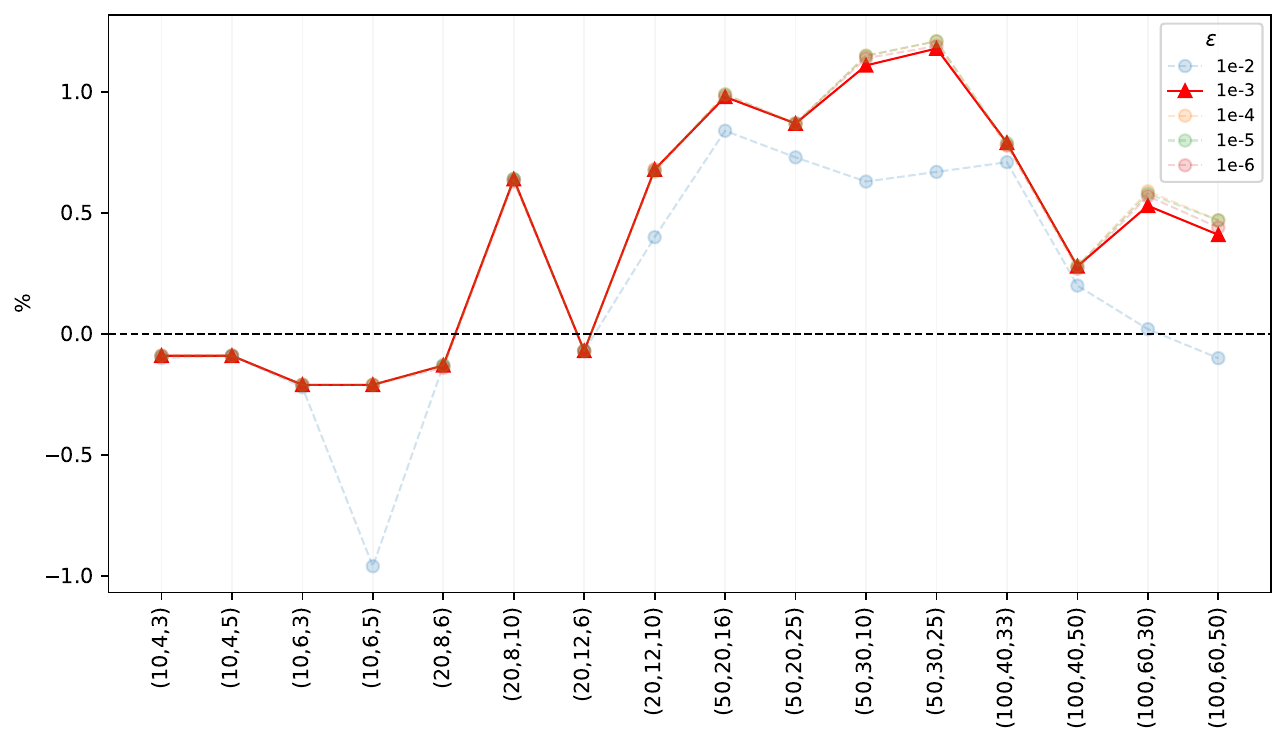} 
        \caption{Average gap of objective value}
        \label{fig:obj_T_2_e}
    \end{subfigure}
    \begin{subfigure}[b]{0.48\textwidth}
        \centering
        \includegraphics[width=\textwidth]{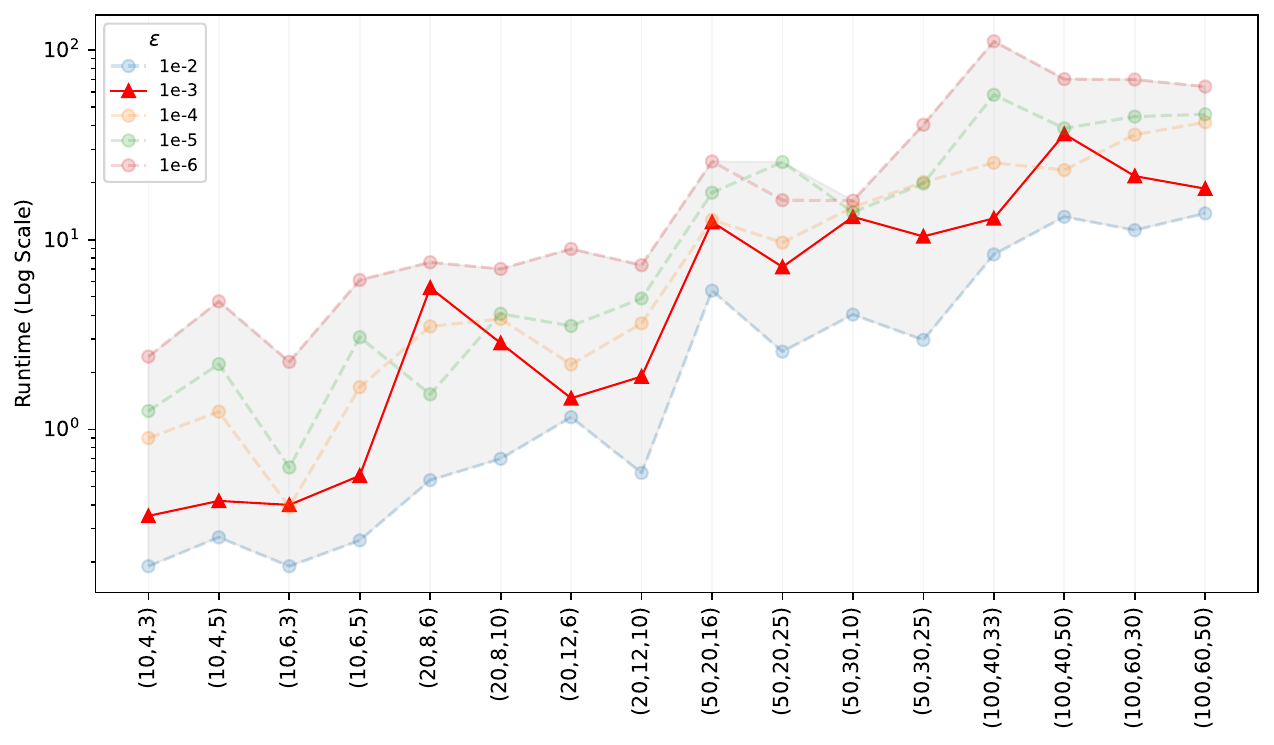} 
        \caption{Average runtime}
        \label{fig:time_T_2_e}
    \end{subfigure}
    
    \caption{Impact of $\epsilon$ to \textbf{LOG-PW} on A\&P instances with $T=2$}
    \label{fig:e_T_2}
\end{figure}

\begin{figure}[!h]
    \centering
    \begin{subfigure}[b]{0.48\textwidth}
        \centering
        \includegraphics[width=\textwidth]{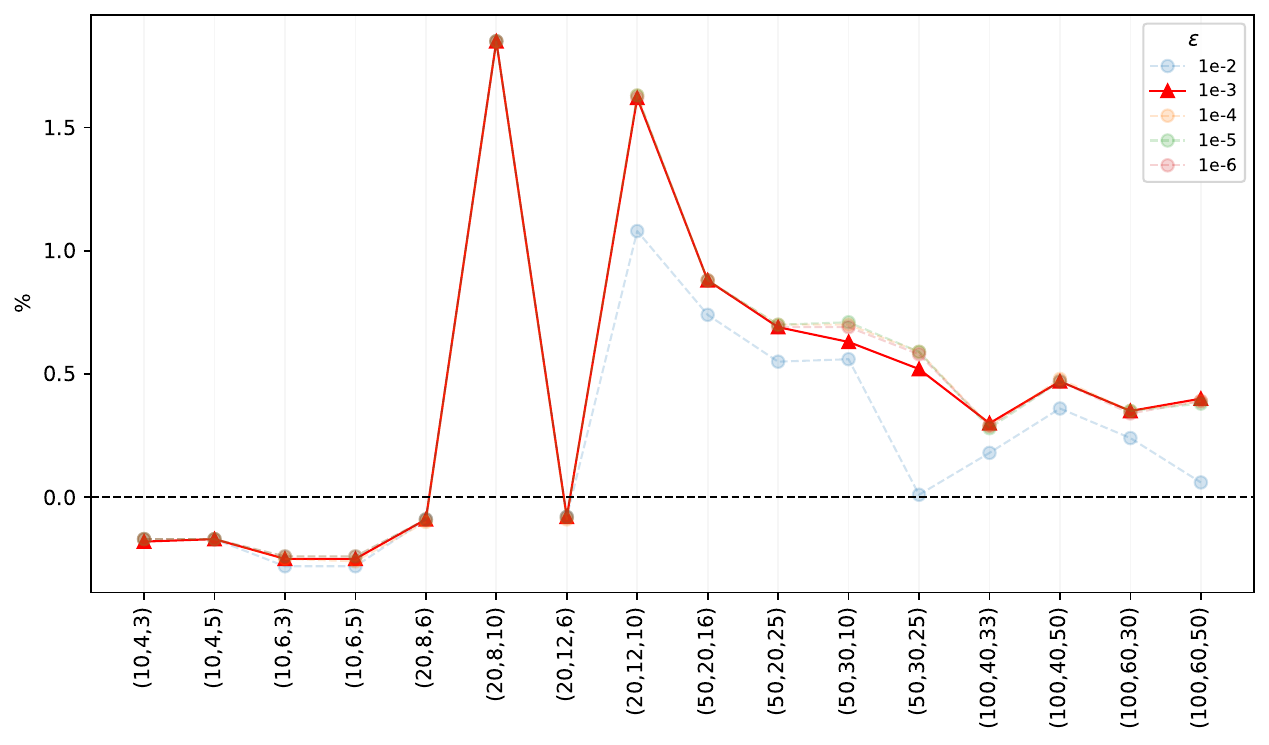} 
        \caption{Average gap of objective value}
        \label{fig:obj_T_5_e}
    \end{subfigure}
    \begin{subfigure}[b]{0.48\textwidth}
        \centering
        \includegraphics[width=\textwidth]{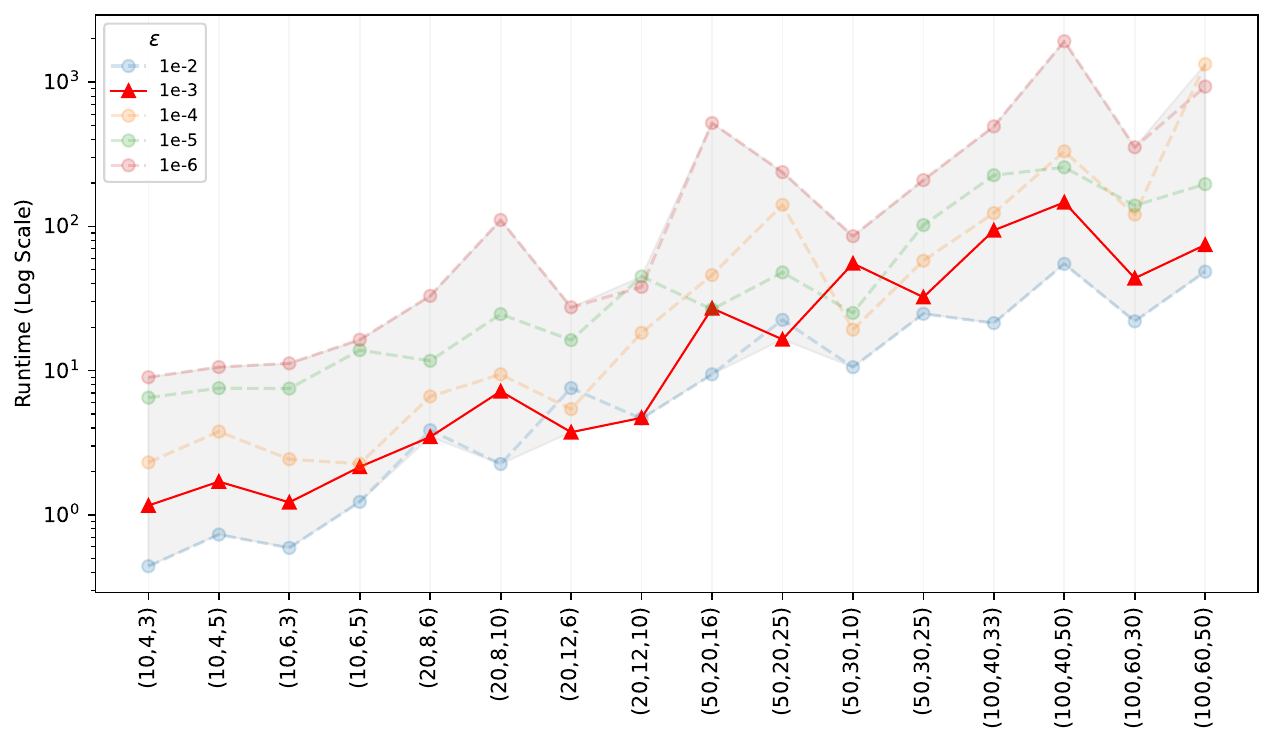} 
        \caption{Average runtime}
        \label{fig:time_T_5_e}
    \end{subfigure}
    
    \caption{Impact of $\epsilon$ to \textbf{LOG-PW} on A\&P instances with $T=5$}
    \label{fig:e_T_5}
\end{figure}

\paragraph{Choice of $K$.} Figures~\ref{fig:K_T_2} and \ref{fig:K_T_5} report the average percentage gap (\%) and average runtime (in log scale) across all A\&P instances, grouped by $(m, C, M)$ for $T=2$ and $T=5$. The results show that for $K \geq 25$, the objective gaps become negligible, suggesting that $K=25$ is a reasonable choice to achieve near-optimal performance across most instances. Based on these findings, we fix $K = 25$ in all subsequent experiments.
\paragraph{Choice of $\epsilon$.} Figures~\ref{fig:e_T_2} and~\ref{fig:e_T_5} illustrate the impact of varying $\epsilon$ on both the solution quality and runtime of our approach (\textbf{LOG-PW}), benchmarked against the baseline provided by \textbf{SCIP}. The objective gap is computed relative to SCIP's output, and runtime is measured in seconds. As observed, reducing $\epsilon$ leads to more accurate approximations (i.e., smaller gaps), but also increases computational cost due to the need for more segments in the PWLA. Interestingly, the results show that for $\epsilon \leq 10^{-3}$, the approximation errors become negligible and the objective gaps converge, yielding nearly identical performance across smaller $\epsilon$ values. This indicates that $\epsilon = 10^{-3}$ is sufficiently small to ensure high-quality solutions without incurring unnecessary computational overhead.

\end{document}